\pdfoutput=1
\RequirePackage{ifpdf}
\ifpdf 
\documentclass[pdftex]{sigma}
\else
\documentclass{sigma}
\fi

\numberwithin{equation}{section}

\newtheorem{Theorem}{Theorem}[section]
\newtheorem*{Theorem*}{Theorem}
\newtheorem{Corollary}[Theorem]{Corollary}
\newtheorem{Lemma}[Theorem]{Lemma}
\newtheorem{Proposition}[Theorem]{Proposition}
 { \theoremstyle{definition}

\newtheorem{Remark}[Theorem]{Remark}
\newtheorem*{Comments}{Comments}
}

\begin{document}
\allowdisplaybreaks

\newcommand{\arXivNumber}{2107.07204}

\renewcommand{\PaperNumber}{068}

\FirstPageHeading

\ShortArticleName{Moduli Spaces for the Fifth Painlev\'e Equation}

\ArticleName{Moduli Spaces for the Fifth Painlev\'e Equation}

\Author{Marius VAN DER PUT and Jaap TOP}
\AuthorNameForHeading{M.~van der Put and J.~Top}
\Address{Bernoulli Institute, Nijenborgh 9, 9747 AG~Groningen, The Netherlands}
\Email{\href{mailto:m.van.der.put@rug.nl}{m.van.der.put@rug.nl}, \href{mailto:j.top@rug.nl}{j.top@rug.nl}}
\URLaddress{\url{https://www.math.rug.nl/~top/}}

\ArticleDates{Received July 15, 2021, in final form September 07, 2023; Published online September 26, 2023}

\Abstract{Isomonodromy for the fifth Painlev\'e equation ${\rm P}_5$ is studied in detail in the context of certain moduli spaces for connections, monodromy, the Riemann--Hilbert morphism, and Okamoto--Painlev\'e spaces. This involves explicit formulas for Stokes matrices and parabolic structures. The rank~4 Lax pair for ${\rm P}_5$, introduced by Noumi--Yamada et al., is shown to be induced by a natural fine moduli space of connections of rank~4. As a~by-product one obtains a polynomial Hamiltonian for ${\rm P}_5$, equivalent to the one of Okamoto.}


\Classification{33E17; 14D20; 14D22; 34M55}

\section{Introduction and summary}

Historically, the list of classical Painlev\'e equations $y''=F(y',y,z)$ was deduced from the property that the only moving singularities
of the solutions are poles. The connection with isomonodromy was quite early recognised. In the papers of Jimbo, Miwa and Ueno \cite{JM,J-M-U, Mi} this theme is developed and moduli spaces for connections and monodromy data are constructed. In particular, isomonodromic families of rank two matrix differential operators $\frac{{\rm d}}{{\rm d}z}+A$ are provided for each of ${\rm P}_1$--${\rm P}_6$. We refer to \cite{C2,C,C3,G,OO1,O1,O,O2} for some theory and
detailed explicit equations. We note that this choice does not do justice to the extensive literature related to Painlev\'e equations.
The family of matrix differential operators for ${\rm P}_3$, ${\rm P}_4$ and for the degenerate fifth Painlev\'e equation ${\rm degP}_5$ has been refined, see \cite{A-P-T, vdP-T,vdP-T3}, to fine moduli spaces of connections on the projective line, which are identified with Okamoto--Painlev\'e spaces. The detailed construction of the moduli spaces supplements and continues some sections of \cite{vdP-Sa}. The moduli spaces for the monodromy data, including Stokes data, are studied in the literature under the name ``wild character varieties''; see, e.g., \cite{Boalch3}.
For the case ${\rm P}_5$, there is a preprint~\cite{PR23} by E.~Paul and J.-P.~Ramis, continuing \cite[Section~3.2]{Pa-Ra}.

The present paper applies a comparable technique for ${\rm P}_5$. This supplements and improves
 \cite[Sections~3.2 and 4.3]{vdP-Sa} involving ${\rm P}_5$. The family of matrix differential operators will be refined to ``natural'' fine moduli spaces $\mathcal{M}(\theta_0,\theta_1,\theta_\infty)$ of connections of rank $2$ on the projective line.

 \subsection[Description of M(theta\_0,theta\_1,theta\_infty) and R\^{}\{geom\}(theta\_0,theta\_1,theta\_infty)]{Description of $\boldsymbol{\mathcal{M}(\theta_0,\theta_1,\theta_\infty)}$ and $\boldsymbol{\mathcal{R}^{\rm geom}(\theta_0,\theta_1,\theta_\infty)}$}\label{sec1.1}

 The objects for the moduli space $\mathcal{M}(\theta_0,\theta_1,\theta_\infty)$
 are connections on a rank two vector bundle of degree $-1$,
 having the (generalized) eigenvalues $\pm \frac{\theta_0}{2}$, $\pm \frac{\theta_1}{2}$ and
 $\frac{tz+\theta_\infty}{2}$, $-\frac{tz+\theta_\infty}{2}-1$ at ${z=0,1}\mbox{ and }\infty$. Further, a variable $u$ with ${\rm e}^{2 \pi {\rm i} u}=t$ is part of the data. Let $\theta :=\theta_\infty +1$.

 There are restrictions on the parameters $\theta_0$, $\theta_1$, $\theta_\infty$ and on the objects of $\mathcal{M}(\theta_0,\theta_1,\theta_\infty)$.
 Moreover, there is an additional structure for $\theta_0\in \mathbb{Z}$ and/or $\theta_1\in\mathbb{Z}$. {\it We give some technical details}.

 For a given connection $(\nabla, \mathcal{V})$ (with $\mathcal{V}$ of degree $-1$), a subobject is a saturated line bundle $\mathcal{L}\subset \mathcal{V}$,
 invariant under $\nabla$. Then $\mathcal{L}$ has local equations $\frac{{\rm d}}{{\rm d}z}+\epsilon_0\frac{\theta_0}{2}$, $\frac{{\rm d}}{{\rm d}z}-\epsilon_1\frac{\theta_1}{2}$
 and $\frac{{\rm d}}{{\rm d}z}-\epsilon_2\frac{tz+\theta}{2}-\frac{1}{2z}$ at the points $z=0,1,\infty$ for certain $\epsilon_0,\epsilon_1,\epsilon_2\in \{\pm 1\}$.
 Then the degree of $\mathcal{L}$ is equal to $-1/2-\epsilon_2\frac{\theta}{2}+\epsilon_1\frac{\theta_1}{2}-\epsilon_0\frac{\theta_0}{2}$.
 We require that any $(\nabla,\mathcal{V})$ has at most one subobject and that the degree of this subobject is $-1$.
 In order to avoid the possibility of subobjects of degree $<-1$, {\it we make in the sequel the following restriction ``${\rm restr}$'' on the parameters $\theta_0$, $\theta_1$, $\theta$:
\[
 -1/2-\epsilon_2\frac{\theta}{2}+\epsilon_1\frac{\theta_1}{2}-\epsilon_0\frac{\theta_0}{2}\in \mathbb{Z}_{<-1}
 \]
is not possible for $(\epsilon_0,\epsilon_1,\epsilon_2)\in \{\pm 1\}^3$.}

For the case $\theta_0\in \mathbb{Z}$, the objects $(\mathcal{V},\nabla)$ are provided with an additional structure called ``parabolic structure''.
 This consist of a one-dimensional $\nabla$-invariant subspace of $\mathcal{V}\otimes \mathbb{C}(\!(z)\!)$. The case $\theta_1\in \mathbb{Z}$ is similar.

\begin{Comments}\quad
\begin{enumerate}\itemsep=0pt
\item[$(1)$] The reason for the choice: ``degree of $\mathcal{V}$ is $-1$'' is the following. In general, $\mathcal{V}\cong O(k)e_1\oplus O(-k-1)e_2$ with $k\geq 0$.
 If $k>0$, then $O(k)e_1$ is a subobject ruled
 out by the restriction on subobjects. This implies $k=0$ and the underlying sheaf
of the connections can be identified with the subsheaf $\mathcal{V}=Oe_1\oplus O(-[\infty])e_2$
of $Oe_1\oplus Oe_2$.

\item[$(2)$] The generic fibre $M$ of a connection $(\mathcal{V},\nabla)$, as above, is a differential module over $\mathbb{C}(z)$ with invariant local
lattices at $z=0,1,\infty$. One can change the local lattices by shifting the $\theta_0$, $\theta_1$, $\theta$ over elements in $2\mathbb{Z}$ and obtain another connection
with generic fibre $M$. This defines in fact a B\"acklund transformation. By doing so, one can obtain parameters for $M$ satisfying ``${\rm restr}$''.

\item[$(3)$] In order to make $\mathcal{M}(\theta_0,\theta_1,\theta_\infty)$ explicit, one computes the objects $(\mathcal{V},\nabla)$ which admit two
subobjects of degree $-1$ (see Section~\ref{4.1}).

\item[$(4)$] Parameters $\{\theta_*\}$ are called {\it exceptional}, (if ``${\rm restr}$'' holds and) if
the equation $-1/2-\epsilon_2\frac{\theta}{2}+\epsilon_1\frac{\theta_1}{2}-\epsilon_0\frac{\theta_0}{2}\ =-1$ has a solution with $\epsilon_2=-1$, $\epsilon_0,\epsilon _1\in \{\pm 1\}$ and a solution
with $\epsilon_2=1$, $\epsilon_0,\epsilon_1\in \{\pm 1\}$. The exceptional parameters are in fact
\begin{gather*}
(\theta=0, -\epsilon_0\theta_0+\epsilon_1\theta_1=-1),\qquad (\theta_0=\pm 1, \theta=\pm \theta_1),
\qquad (\theta_1=\pm 1, \theta=\pm \theta_0).
\end{gather*}
In order to define the Riemann--Hilbert morphism for exceptional parameters, one refines the space of connection.
Let $\mathcal{M}(\theta_0,\theta_1,\theta_\infty, \epsilon_2=-1)$ denote the locus in $\mathcal{M}(\theta_0,\theta_1,\theta_\infty)$
where the reducible objects correspond to $\epsilon_2=-1$. Further, $\mathcal{M}(\theta_0,\theta_1,\theta_\infty, \epsilon_2=1)$ is
the locus where the reducible objects correspond to $\epsilon_2=1$.
A moduli space $\mathcal{R}^{\rm geom}(\theta_0,\theta_1,\theta_\infty)$ for the monodromy
is constructed from the identity ${\rm mon}_0\cdot {\rm mon}_1={\rm mon}_\infty$, where the $\{{\rm mon}_*\}$ are the local topological
monodromies at~$0$,~$1$,~$\infty$, and from the monodromy identity which expresses ${\rm mon}_\infty$ as a product of the formal monodromy and two Stokes matrices. This moduli space depends in a delicate way on the possibilities for subobjects
of the space of connections. In case of exceptional parameters, one defines two spaces
 $\mathcal{R}^{\rm geom}(\theta_0,\theta_1,\theta_\infty, \epsilon_2=-1)$ and $\mathcal{R}^{\rm geom}(\theta_0,\theta_1,\theta_\infty, \epsilon_2=1)$
which correspond to the similar refinement for the space of connections.
 The above monodromy space is also provided with parabolic structures for the cases $\theta_0\in \mathbb{Z}$ and/or $\theta_1\in \mathbb{Z}$.
 This structure consists of an invariant line for the local monodromy at $z=0$ and/or the
 local monodromy at $z=1$.
\end{enumerate}
\end{Comments}

 In \cite[Section~3.2]{Pa-Ra}, the moduli space $\mathcal{R}^{\rm geom}(\theta_0,\theta_1,\theta_\infty)$ is seen as a representation of the wild fundamental groupoid and the monodromy identity is also present in that context.

The main result (Theorem~$\ref{2.5}$):
{\it If the parameters satisfy ``${\rm restr}$'' and are not exceptional, then the extended Riemann--Hilbert morphism
\[{\rm RH}^+\colon \ \mathcal{M}(\theta_0,\theta_1,\theta_\infty)\rightarrow
\mathcal{R}^{\rm geom}(\theta_0,\theta_1,\theta_\infty)\times \mathbb{C}\]
is an isomorphism. Moreover, for exceptional parameters the same holds after adding $\epsilon_2=-1$ or $\epsilon_2=1$ to both spaces.}

This implies the Painlev\'e property for the solutions of ${\rm P}_5$ and provides a
comparison of $\mathcal{M}(\theta_0,\theta_1,\theta_\infty)$ with the Okamoto--Painlev\'e spaces.

Some ``natural'' moduli spaces of connections of rank 4, related to ${\rm P}_5$, are introduced. The corresponding Lax pairs are identified with the Lax pairs of Noumi--Yamada for ${\rm P}_5$.
All B\"acklund transformations for ${\rm P}_5$ are explicit from these Lax pairs.

{\it Now we discuss some more details of the paper}. The construction of moduli spaces starts in Section~\ref{Section1} by defining a set
${\bf S}(\theta_0,\theta_1,\theta_\infty)$ of pairs $(M,u)$, where $M$ is a differential module of dimension 2 over $\mathbb{C}(z)$ and
$u \in \mathbb{C}$. Further, $M$ has regular singular points $z=0$ and $z=1$ with local eigenvalues $\pm \frac{\theta_0}{2}$ and $\pm \frac{\theta_1}{2}$. Moreover, the point $z=\infty$ is an irregular singularity with (generalized) eigenvalues $\pm \frac{tz+\theta_\infty}{2}$ and $t={\rm e}^{2\pi {\rm i} u}$.

The choice of the parameters $\{\theta_*\}$ and the parabolic structures (in case $\theta_0$ and/or $\theta_1$ are
in~$\mathbb{Z}$) determines a connection with generic fibre $M$. We assume again that $\{\theta_*\}$ satisfies ``${\rm restr}$''. Further, we allow for at most
one proper submodule $L\subset M$. If $L$ exists, then the degree of the corresponding saturated line bundle $\mathcal{L}$ is required to be
$-1$. In case of exceptional parameters we refine this set into the two sets ${\bf S}(\theta_0,\theta_1,\theta_\infty, \epsilon_2=1)$ and
${\bf S}(\theta_0,\theta_1,\theta_\infty, \epsilon_2=-1)$.

The monodromy data for objects of ${\bf S}(\theta_0,\theta_1,\theta_\infty)$ are given by local monodromy matrices, Stokes matrices and parabolic structure. In Section~\ref{Section2}, the
construction of the fine moduli space $\mathcal{R}^{\rm geom}(\theta_0,\theta_1,\theta_\infty)$
is explained. In the presence of reducible objects this space required a more refined description, see Section~\ref{sec1.1} and page~\pageref{twee}.

It is shown that $\mathcal{R}^{\rm geom}(\theta_0,\theta_1,\theta_\infty)$ is simply connected. It
is a resolution of the monodromy space obtained in \cite[Sections~3.2 and~4.3]{vdP-Sa}, if one assumes that the parameters
satisfy ``${\rm restr}$'' and are not exceptional; indeed, see Remark~\ref{remark3.3}\,(b).

In the exceptional case, the monodromy space is refined into spaces $\mathcal{R}^{\rm geom}(\theta_0,\theta_1,\theta_\infty, \epsilon_2=1)$ and $\mathcal{R}^{\rm geom}(\theta_0,\theta_1,\theta_\infty, \epsilon_2=-1)$.

 As always we assume ``${\rm restr}$''. If the parameters are not exceptional, then the natural map
${\bf S}(\theta_0,\theta_1,\theta_\infty)\rightarrow \mathcal{R}^{\rm geom}(\theta_0,\theta_1,\theta_\infty)(\mathbb{C})\times \mathbb{C}$ is a bijection.
 If the parameters are exceptional, then a similar natural bijection is valid after adding $\epsilon_2=1$ or~$\epsilon_2=-1$ to both sides.

In Section~\ref{Section3}, one considers the vector bundle $\mathcal{V}=Oe_1\oplus O(-[\infty])e_2$,
subbundle (of degree~$-1$) of the free vector bundle $Oe_1\oplus Oe_2$ on $\mathbb{P}^1$.
The data for every object $(M,u)\in {\bf S}(\theta_0,\theta_1,\theta_\infty)$ produce a unique connection
$\nabla(M)\colon \mathcal{V}\rightarrow \Omega([0]+[1]+2[\infty] )\otimes \mathcal{V}$ with prescribed data at the singular points $0$, $1$, $\infty$. The generic fibre of $\nabla(M)$ is the differential module $M$ (including data).
This leads to a moduli functor and a fine moduli space $\mathcal{M}(\theta_0,\theta_1,\theta_\infty)$.

 It is shown that the extended Riemann--Hilbert morphism ${\rm RH}^+$ is an isomorphism if
 ${\theta_0,\theta_1 \! \not \in \! \mathbb{Z}}$, the parameters satisfy ``${\rm restr}$'' and are not exceptional. In the exceptional case, the statement
 remains valid after adding $\epsilon_2=-1$ or $ \epsilon_2=1$ to all spaces.

 The reducible locus, corresponding to the reducible $M\in {\bf S}(\theta_0,\theta_1,\theta_\infty)$ is computed
 in Section~\ref{Section4}. In the presence of reducible objects the definition of $\mathcal{M}$ has to be refined, see Section~\ref{sec1.1} and page~\pageref{Mrefined} for details. Isomonodromy in the sublocus of $\mathcal{M}(\theta_0,\theta_1,\theta_\infty)$
corresponding to reducible modules produces the
 Riccati solutions for ${\rm P}_5$. The parabolic structure is again studied in Section~\ref{Section5} in order to prove the main result,
 Theorem~\ref{2.5}, namely:
\textit{The extended Riemann--Hilbert morphism is an isomorphism for parameters satisfying ``${\rm restr}$'' $($and adding
 $\epsilon_2=-1$ or $\epsilon_2=1$ for the exceptional cases$)$.}

In Section~\ref{Section6}, a standard computation with Lax Pairs produces the usual differential equation for~${\rm P}_5$.
From the extended Riemann--Hilbert map ${\rm RH}^+$ one can read off and find the formulas for some B\"acklund transformations.
 However, some are missing. Noumi and Yamada constructed a Lax pair for ${\rm P}_5$ from which one can read off all B\"acklund transformations.

 Our contribution in Section~\ref{Section7} to this is to provide a natural moduli space of connections
 $\mathcal{M}_4$ and a corresponding monodromy space $\mathcal{R}_4$, which produce
the Lax pair of Noumi--Yamada. The moduli space $\mathcal{M}_4$ is associated to the set of differential modules $M$ of dimension 4 over $\mathbb{C}(z)$ such that:
\begin{itemize}\itemsep=0pt
\item[$({\rm a})$] $\Lambda ^4M$ is trivial,
\item[$({\rm b})$] $z=0$ is regular singular,
\item[$({\rm c})$] $z=\infty$ is irregular singular, and its generalized eigenvalues are $z^{1/2}+tz^{1/4}$ and its conjugates.
\end{itemize}

 A direct computation of a matrix differential operator for $\mathcal{M}_4$ seems rather difficult.
Instead, one gives a module $M\in \mathcal{M}_4$ the interpretation of a
 differential module $N$ over $\mathbb{C}\big(z^{1/4}\big)$ provided with a symmetry $\sigma$ of order~4.
 Explicit formulas for $N$ and $M:=N^{\langle\sigma\rangle}$ are inspired by \cite[Section~12.5]{vdP-S}. The resulting Lax pair is identical to the one in \cite{S-H-C}. A calculation of the Lax pair for $N$ instead of $M$ produces at once a polynomial Hamiltonian for ${\rm P}_5$ which is, up to a~linear change of coordinates, also present in Okamoto's work.

 \subsection*{Comparison with work of Ph.~Boalch}
In the remainder of this introduction, we follow a suggestion by the referee of earlier versions of this text, by presenting details of the subtle comparison of our work with the paper \cite{Boalch} by Ph.~Boalch.
We note that the paper \cite{Boalch2} extends the constructions and results of \cite{Boalch} to connections on algebraic curves of higher genus. We will not discuss~\cite{Boalch2}.

Firstly, we describe some details of \cite{Boalch}, relevant for the comparison. The data $\mathbb{A}$ denotes a~set of singular points $\{a_1,\dots ,a_m\}\subset \mathbb{P}^1$ together with for each point the type of the singularity. For a local coordinate $w$ at a singular point $a_j$, the differential operator is supposed to have the form $\frac{{\rm d}}{{\rm d}w}+\frac{A_k}{w^k}+\frac{A_{k-1}}{w^{k-1}}+\cdots$ with $k\geq 1$ such that $A_k$ has distinct eigenvalues for $k\geq 2$ and the eigenvalues are also distinct modulo $\mathbb{Z}$ in case $k=1$. The ``type'' of the singularity $a_j$ is the principle part of this operator.

 Associated to $\mathbb{A}$ are moduli sets $\mathcal{M}^*(\mathbb{A})\subset \mathcal{M}(\mathbb{A})$ and $M(\mathbb{A})$.
 The first one is the set of isomorphy classes of connections
with data $\mathbb{A}$ on a trivial vector bundle of rank 2 on $\mathbb{P}^1$. The second moduli set is the same but now with ``trivial vector bundle'' replaced
by ``vector bundle of degree zero''. Finally, $M(\mathbb{A})$ is a moduli set for the analytic data, i.e., monodromy, Stokes matrices (and links)
associated to $\mathbb{A}$.

The set $\mathcal{M}^*(\mathbb{A})$ is identified with $O_1\times \cdots \times O_m$ modulo the action by conjugation of ${G={\rm GL}_n}$. Here each
$O_j$ is a $G_{k_j}:={\rm GL}_n\big(\mathbb{C}[\zeta]/\zeta^{k_j}\big)$-orbit in the Lie algebra of the group $G_{k_j}$. Since (co)adjoint orbits have a natural symplectic structure, the same holds for $O_1\times \cdots \times O_m$.
Now~$\mathcal{M}^*(\mathbb{A})$, regarded as variety, is defined as the categorical quotient $(O_1\times \cdots \times O_m)//G$.
We note in passing that this categorical quotient need not be a geometric quotient and that $(O_1\times \cdots \times O_m)//G$ can have singularities.
A priori, the moduli set $\mathcal{M}(\mathbb{A})$ is not given a~structure of algebraic variety.

The $C^\infty$-method of B.~Malgrange et al.\ to deal with
the Stokes phenomenon provides the structure of an algebraic variety on $M(\mathbb{A})$ and its symplectic structure. We omit details.

A main result of \cite{Boalch} states that the Riemann--Hilbert morphism $v\colon \mathcal{M}^*(\mathbb{A})\rightarrow M(\mathbb{A})$
is a~symplectic, analytic isomorphism onto its image, which is a dense open subset of $M(\mathbb{A})$. Moreover, $v$ extends to a {\it bijection} $v^{\rm ext}\colon \mathcal{M}(\mathbb{A})\rightarrow M(\mathbb{A})$.
The bijection $v^{\rm ext}$ provides $\mathcal{M}(\mathbb{A})$ with the structure of
symplectic analytic variety. The restriction of these structures to
$\mathcal{M}^*(\mathbb{A})$ coincides with the already constructed data.

A further step is to replace the fixed data $\mathbb{A}$ by deformation data, say consisting of a ``time variable'' $t\in T$.
 Then the set $M(\mathbb{A})$ is closely related to our construction of $\mathcal{R}\times T$, where~$\mathcal{R}$ is ``our'' monodromy space:
this is the transparent structure of~$\mathcal{R}$ derived by
the technique of J.~Martinet, J.-R.~Ramis, B.L.J.~Braaksma et al.\ of (multi)summation applied to Stokes theory. We believe that it can be verified that
the algebraic structures and the symplectic structures of~$M(\mathbb{A})$ and $\mathcal{R}\times T$ coincide.

{\it Consider now the case ${\rm P}_5$} (not made explicit in \cite{Boalch}).
 Then $\mathbb{A}$ corresponds to singularities~$0$,~$1$,~$\infty$ with eigenvalues
 $\pm \frac{\theta_0}{2}$, $\pm \frac{\theta_1}{2}$, $\pm \frac{tz+\theta_\infty}{2}$ with $\theta_0, \theta_1\not \in \mathbb{Z}$. The connections on the trivial bundle
 produce a variety $\mathcal{M}^*(\mathbb{A})$ as a quotient $O_0\times O_1\times O_\infty //G$, where $G={\rm GL}_2$.

 If this is seen as categorical quotient, then the space has singularities, due to reducibility, for special values of the $\theta$'s. In such a case it is not a geometric quotient. After leaving out the
 locus corresponding to reducible modules, there is a good geometric quotient. We assume that this is a good interpretation of
Boalch's moduli space $\mathcal{M}^*(\mathbb{A})$.

Not every differential module $M$ with the above data can be realized as connection on a~trivial bundle.
However, every $M$ can be represented by a connection on a vector bundle with degree~$-1$.
After putting a restriction on the parameters and the reducible modules, ``our'' space $\mathcal{M}(\theta_0,\theta_1,\theta_\infty)$ is obtained as good geometric quotient under the group $G$. Moreover, a~natural
 symplectic structure can be deduced. We conclude that
$\mathcal{M}^*(\mathbb{A})$ is a Zariski open subset of $\mathcal{M}(\theta_0,\theta_1,\theta_\infty)$. It can probably be identified with the complement of the tau-divisor; see \cite[Section~3.3.4]{vdP-T3} for
details in analogous cases.

We note that for all cases ${\rm P}_1$--${\rm P}_6$ the choice of a vector bundle of rank 2 and degree $-1$ leads to a good moduli space of connections.

We observe that every $M$, as above, can be realized as connection on a vector bundle of degree 0. Thus $\mathcal{M}(\mathbb{A})$
 is closely related to $\mathcal{M}(\theta_0,\theta_1,\theta_\infty)$.

 \section[The set S(theta\_0,theta\_1,theta\_infty) of differential modules]{The set $\boldsymbol{{\mathbf S}(\theta_0,\theta_1,\theta_\infty)}$ of differential modules}\label{Section1}
 Fix parameters $\theta_0,\theta_1,\theta_\infty \in \mathbb{C}$. Put $\theta=\theta_\infty +1$ and $t={\rm e}^{2\pi {\rm i} u}$.
 We consider a set $\mathbf{S}(\theta_0,\theta_1,\theta_\infty)$ of pairs $(M,u)$ (up to isomorphism) with $u\in \mathbb{C}$ and $M$ a differential module over $\mathbb{C}(z)$, defined by the following properties:
\begin{enumerate}\itemsep=0pt
\item[$(1)$] $M$ has dimension two; $\det M=\Lambda^2M$ is the trivial one-dimensional differential module; $M$~has at most singularities at $z=0,1,\infty$.

\item[$(2)$] The point $z=0$ is regular singular and $\mathbb{C}(\!(z)\!)\otimes M$
 can be represented by the differential operator
 $\frac{{\rm d}}{{\rm d}z}+\frac{1}{z}A_0$, where $A_0$ is a $(2\times 2)$-matrix with entries in $\mathbb{C}[\![z]\!]$ and
 $A_0 \bmod (z)$ has eigenvalues $\pm \frac{\theta_0}{2}$.

\item[$(3)$] The point $z=1$ is regular singular and $\mathbb{C}(\!(z-1)\!)\otimes M$ can be represented by the operator
 $\frac{{\rm d}}{{\rm d}z}+\frac{1}{z-1}A_1$, where $A_1$ is a $(2\times 2)$-matrix with entries in $\mathbb{C}[\![z-1]\!]$
 and $A_1\bmod(z-1)$ has eigenvalues $\pm \frac{\theta_1}{2}$.

\item[$(4)$] The point $z=\infty$ is irregular singular and $\mathbb{C}\big(\!\big(z^{-1}\big)\!\big)\otimes M$ can be represented by the operator
\[\frac{{\rm d}}{{\rm d}z}+ \frac{1}{z}\cdot \left(\begin{matrix} \dfrac{tz+\theta_\infty}{2}& 0\\ 0& -\dfrac{tz+\theta_\infty}{2} \end{matrix}\right).\]
\end{enumerate}

 \begin{Comments} The ``datum'' $u$ is attached to the module $M$ in order to distinguish the two
 generalized local exponents $\pm \frac{tz+\theta_\infty}{2}$ at $z=\infty$. The aim is to
 construct a moduli space $\mathcal{M}(\theta_0,\theta_1,\theta_\infty)$ of connections
 with a natural bijection $\mathcal{M}(\theta_0,\theta_1,\theta_\infty)(\mathbb{C})\rightarrow {\bf S}(\theta_0,\theta_1,\theta_\infty)$.

 For the monodromy, we will construct a moduli space $\mathcal{R}^{\rm geom}(\theta_0,\theta_1,\theta_\infty)$.
 The term $u$, and not only $t$, is needed for the definition of the Riemann--Hilbert morphism
 ${\rm RH}\colon \mathcal{M}(\theta_0,\theta_1,\theta_\infty)\rightarrow
 \mathcal{R}^{\rm geom}(\theta_0,\theta_1,\theta_\infty)$.
\end{Comments}
\begin{enumerate}\itemsep=0pt
\item[$(5)$] The following technical condition on reducible modules $M$:
A reducible module $M$, satisfying (1)--(4), belongs to ${\bf S}(\theta_0,\theta_1,\theta_\infty)$ if it has only one
submodule $L$ of dimension 1 and moreover $L$ satisfies a certain condition. Namely, $L$
gives rise to a tuple $(\epsilon_0,\epsilon_1,\epsilon_2)\in \{\pm 1\}^3$
and an integer $k=-\frac{1}{2}-\epsilon_2\frac{\theta}{2}+\epsilon_1\frac{\theta_1}{2}-\epsilon_0\frac{\theta_0}{2}$.
 We require that $k=-1$. Moreover, we avoid the possibility of an integer $k<-1$. This leads to the assumption of
condition ``${\rm restr}$'' (see Section~\ref{sec1.1}, Corollary~\ref{2.*} and Remark~\ref{1.3} for details).
Parameters $\{\theta_*\}$ are {\it exceptional} if both $\epsilon_2=-1$ and $\epsilon_2=1$ occurs for reducible objects with degree $-1$.
For exceptional parameters, the definition of ${\bf S}(\theta_0,\theta_1,\theta_\infty)$ is refined by
introducing the sets ${\bf S}(\theta_0,\theta_1,\theta_\infty, \epsilon_2=-1)$
and ${\bf S}(\theta_0,\theta_1,\theta_\infty, \epsilon_2=1)$. For the first one, we require additionally for reducible modules that $\epsilon_2=-1$.
For the second set, we require additionally that $\epsilon_2=1$.
{\it For not exceptional parameters, we keep the notation ${\bf S}(\theta_0,\theta_1,\theta_\infty)$.}
\end{enumerate}
 \begin{Comments}
For the above submodule $L$, the local data at the points $0$, $1$, $\infty$ are $\frac{{\rm d}}{{\rm d}z}+\epsilon_0\frac{\theta_0}{2z}$,
 $\frac{{\rm d}}{{\rm d}z}-\epsilon_1\frac{\theta_1}{2(z-1)}$ and $\frac{{\rm d}}{{\rm d}z}-\epsilon_2\frac{tz+\theta}{2z} -\frac{1}{2z}$.
Furthermore, $k$ is the degree of the line bundle $\mathcal{L}$ corresponding to $L$ and the lattices described by the local data.

 Split modules $M$ (i.e., direct sums of submodules of dimension 1), satisfying (1)--(4), may produce a bad singularity or
non-separated moduli spaces. Therefore, split modules are not admitted in any set ${\bf S}(\theta_0,\theta_1,\theta_\infty)$.

 If one admits (for given exceptional parameters) reducible modules for both cases $\epsilon_2=-1$ and $\epsilon_2=1$,
then the construction of the moduli space for the connections and/or the one for the monodromy could produce non-separated spaces.
The exceptional parameters are
\[ (\theta=0, -\epsilon_0\theta_0+\epsilon_1\theta_1=-1),\qquad (\theta_0=\pm 1, \theta=\pm \theta_1),
\qquad (\theta_1=\pm 1, \theta=\pm \theta_0), \]
see Sections~\ref{sec1.1}, \ref{Section2}\,(v)\,(c) and~\ref{4.1} for more details.
\end{Comments}

 For $\theta_0\in \mathbb{Z}$ and/or $\theta_1\in \mathbb{Z}$,
 another refinement of the definition of ${\bf S}(\theta_0,\theta_1,\theta_\infty)$ is needed. This consists of the additional
 data of one or two lines, namely:
\begin{enumerate}\itemsep=0pt
\item[$(6)$] For $\theta _0\in \mathbb{Z}$, it is a line
 $\mathbb{C}v\subset \mathbb{C}(\!(z)\!)\otimes M$ such that $D(v)=-\frac{|\theta_0|}{2} v$.
\item[$(7)$] For $\theta_1\in \mathbb{Z}$, it is a line
 $\mathbb{C}v\subset \mathbb{C}(\!(z-1)\!)\otimes M$ such that $D(v)=-\frac{|\theta_1|}{2} v$.
The lines are called {\it eigenlines} and the additional data is called a {\it parabolic structure} or a {\it level structure}.
\end{enumerate}

{\it Explanation of the conditions} (6) {\it and }(7).
Consider a differential module $M$ over $\mathbb{C}(z)$ satisfying (1)--(5). Write $D\colon M\rightarrow M$
for its differential operator with $D(f\cdot m)=\frac{{\rm d}f}{{\rm d}z}\cdot m+f\cdot D(m)$. We want to attach to $M$ an
``invariant lattice at $z=0$ with data $\theta_0$''. This is a $\mathbb{C}[\![z]\!]$-lattice
$\Lambda\subset \mathbb{C}(\!(z)\!)\otimes M$ such that $zD(\Lambda)\subset \Lambda$ and
$zD$ induces on $\Lambda/z\Lambda$ a linear operator with eigenvalues $\pm \frac{\theta_0}{2}$. In general $\Lambda$
is unique, except possibly for $\theta_0\in \mathbb{Z}_{\neq 0}$. In the latter case, the parabolic structure, namely
the eigenline
 $\mathbb{C}v$, leads to a unique possibility for $\Lambda$. Without the parabolic structure the (to be constructed) moduli space of connections would have singularities or its set of closed points would differ from
 ${\bf S}(\theta_0,\theta_1,\theta_\infty)$. The parabolic structure produces a natural resolution of the singularities that one would have otherwise. The technical Lemma~\ref{2.1} provides the verification of the above statements. In Section~\ref{Section5}
there is another formulation of this lemma.

\begin{Lemma}\label{2.1} Consider differential modules $M=(M,D)$ over $\mathbb{C}(\!(z)\!)$ with
$D(f\cdot m)=z\frac{{\rm d}f}{{\rm d}z}\cdot m+f\cdot D(m)$. Let $L=\mathbb{C}(\!(z)\!)e$ be defined by $De=\frac{1}{2}e$.
Let $N$ be a differential module with $\dim N=2$, $\det N$ is trivial and $N$ has a $D$-invariant
lattice $\Lambda$ such that $D$ has eigenvalues $\pm \frac{\theta}{2}$ for its action on $\Lambda/z\Lambda$.
\begin{enumerate}\itemsep=0pt
\item[$(1)$] The $D$-invariant lattice $\Lambda$ is unique except in the cases:
\begin{enumerate}\itemsep=0pt
\item[$({a})$] $\theta\in \mathbb{Z}_{\neq 0}$ is even and $N$ is trivial,
\item[$({b})$] $\theta\in \mathbb{Z}_{\neq 0}$ is odd and $L\otimes N$ is trivial.
\end{enumerate}
\item[$(2)$] Let $\theta =m\in \mathbb{Z}_{>0}$ and suppose that {\rm (1a)} or {\rm (1b)} holds.
Let $\mathcal{LA}$ denote the set of $D$-invariant lattices $\Lambda$ such that $D$ has eigenvalues
$\pm \frac{m}{2}$ on $\Lambda/z\Lambda$. Let $\mathcal{LI}$ denote the set of lines $\mathbb{C}v\neq 0$ in $N$
such that $D(v)=-\frac{m}{2}v$. Every lattice $\Lambda \in \mathcal{LA}$ contains a unique
$\mathbb{C}v\in \mathcal{LI}$. This defines a bijection $\mathcal{LA}\rightarrow \mathcal{LI}$.
Fix a basis $b_1$, $b_2$ of $N$ such that $D(b_1)=-\frac{m}{2}b_1,\ D(b_2)=\frac{m}{2}b_2$. The set
$\mathcal{LI}$ identifies
with the projective complex line $\mathbb{P}(\mathbb{C} b_1+\mathbb{C}z^{-m}b_2)$.
\item[$(3)$]Suppose $\theta \in \mathbb{Z}_{\geq 0}$. There is a unique line $\mathbb{C}v\neq 0$ in $N$ with
$D(v)=-\frac{\theta}{2}v$, expect for cases~$(1a)$, $(1b)$ and the case where $N$ is trivial and $\theta=0$.
In the last case the set of those lines forms the projective line $\mathbb{P}(\{n\in N| D(n)=0\})$.
\end{enumerate}
\end{Lemma}

\begin{proof} {\it Suppose $\theta \neq 0$}. The choice of a $D$-invariant lattice for $N$ produces a differential operator
$z\frac{{\rm d}}{{\rm d}z}+A$ where $A$ has entries in $\mathbb{C}[\![z]\!]$ and $A\equiv \big(\begin{smallmatrix} \lambda & 0\\ 0& -\lambda\end{smallmatrix}\big) \bmod (z)$
and $\lambda =-\frac{\theta}{2}$. By conjugation with a~matrix $U=1+U_1z+U_2z^2+\cdots $ (all $U_i$ are
$(2\times 2)$-matrices, $U_0:=1$) we try to obtain $z\frac{{\rm d}}{{\rm d}z}+\big(\begin{smallmatrix} \lambda& 0\\ 0 & -\lambda\end{smallmatrix}\big)$. This leads to
$U\big(z\frac{{\rm d}}{{\rm d}z}+A\big)=\big(z\frac{{\rm d}}{{\rm d}z}+ \big(\begin{smallmatrix}\lambda& 0\\ 0& -\lambda\end{smallmatrix}\big)\big)U$. Write $A=A_0+A_1z+A_2z^2+\cdots$ with
$A_0=\big(\begin{smallmatrix}\lambda& 0\\ 0& -\lambda\end{smallmatrix}\big)$ and all $A_j$ are $(2\times 2)$-matrices. Then one obtains equations
\[
\sum _{n\geq 1}\sum_{i=1}^n(U_{n-i}A_i) z^n =\sum _{n\geq 1} (nU_n+A_0U_n-U_nA_0)z^n
\]
 for the matrices $U_i$.
The map $U_n\mapsto nU_n+A_0U_n-U_nA_0$ sends $\big(\begin{smallmatrix} x_1& x_2\\ x_3& x_4\end{smallmatrix}\big)$
to $\big(\begin{smallmatrix} nx_1& (n+2\lambda)x_2 \\ (n-2\lambda) x_3 & nx_4 \end{smallmatrix}\big)$.

 {\it If $\theta =-2\lambda \not \in \mathbb{Z}$}, there is a unique solution $U=1+U_1z+U_2z^2+\cdots$.
 The given lattice has a~basis $b_1$, $b_2$ with $D(b_1)=\lambda b_1$, $D(b_2)=-\lambda b_2$.

 Any $D$-invariant
 lattice $\tilde{\Lambda}$ such that $D$ has eigenvalues $\pm \lambda$ on $\tilde{\Lambda}/z\tilde{\Lambda}$
 has a basis $\tilde{b}_1$, $\tilde{b}_2$ with $D\big(\tilde{b}_1\big)=\lambda \tilde{b}_1$, $D\big(\tilde{b}_2\big)=-\lambda \tilde{b}_2$.
 The matrix $\big(\begin{smallmatrix}x_1& x_2\\ x_3& x_4 \end{smallmatrix}\big)$
 relating the two bases has the property
 \[\left(\begin{matrix} x_1& x_2\\ x_3& x_4\end{matrix}\right)\left(z\frac{{\rm d}}{{\rm d}z}+
\left(\begin{matrix}\lambda& 0\\ 0& -\lambda \end{matrix}\right) \right)
 =\left(z\frac{{\rm d}}{{\rm d}z}+\left(\begin{matrix} \lambda & 0\\ 0& -\lambda\end{matrix}\right) \right)\left(\begin{matrix} x_1 & x_2\\ x_3& x_4\end{matrix}\right) .\]
 Using $\big(\begin{smallmatrix} x_1& x_2\\ x_3& x_4\end{smallmatrix}\big) \in {\rm GL}_2(\, \mathbb{C}(\!(z)\!)\, )$ one obtains $x_1,x_4\in \mathbb{C}^*$ and
 $x_2=x_3=0$. Thus $\tilde{\Lambda}$ coincides with the given $D$-invariant lattice.

{\it Case $\theta=-2\lambda \in \mathbb{Z}_{\neq 0}$}. We suppose for convenience that $\theta =m\in \mathbb{Z}_{>0}$.
The above system of equations shows that the equations for the $U_n$ have unique solutions for $n\neq m$.
There might not be a solution for $U_m$. This implies that the given lattice has a basis $b_1$, $b_2$ such that $D(b_1)=\lambda b_1$ and $D(b_2)=-\lambda b_2+az^{m}b_1$ where $a\in \mathbb{C}$.

{\it The subcase $a\neq 0$}. We may normalize to $a=1$.
A straightforward calculation shows that there is no element $v\in N-\mathbb{C}z^mb_1$ with $D(v)=-\lambda v$. Any $D$-invariant lattice $\tilde{\Lambda}$ in $N$ such that~$D$ has eigenvalues $\pm \lambda$
on $\tilde{\Lambda}/z\tilde{\Lambda}$ has a basis $\tilde{b}_1$, $\tilde{b}_2$ with similar properties. Comparing
the two lattices produces (as before) an equality
\[
\left(\begin{matrix} x_1& x_2\\ x_3& x_4 \end{matrix}\right) \left(z\frac{{\rm d}}{{\rm d}z}+
\left(\begin{matrix} \lambda& z^m\\ 0& -\lambda\end{matrix}\right) \right)
 =\left(z\frac{{\rm d}}{{\rm d}z}+\left(\begin{matrix}\lambda& z^m\\ 0& -\lambda\end{matrix}\right) \right)\left(\begin{matrix} x_1& x_2\\ x_3& x_4\end{matrix}\right)
 \]
 with
 $\big(\begin{smallmatrix} x_1& x_2\\ x_3& x_4 \end{smallmatrix}\big) \in {\rm GL}_2(\mathbb{C}(\!(z)\!) )$. This leads to the matrix equation
 \[
z\frac{{\rm d}}{{\rm d}z} \left(\begin{matrix} x_1& x_2\\ x_3& x_4\end{matrix}\right) =
 -m\cdot \left[\left(\begin{matrix} x_1& x_2\\ x_3& x_4\end{matrix}\right) ,\left(\begin{matrix} \frac{1}{2}& 0\\ 0& -\frac{1}{2}\end{matrix}\right) \right]+
 z^m\cdot \left[\left(\begin{matrix} x_1& x_2\\ x_3& x_4\end{matrix}\right) ,\left(\begin{matrix} 0& 1\\ 0& 0\end{matrix}\right) \right].\]
The solutions of this system of equations are $x_1=x_4\in \mathbb{C}^*$, $x_3=0$ and $x_2=c\cdot z^m$
with $c\in \mathbb{C}$. Hence
$\tilde{\Lambda}$ coincides with the given $D$-invariant lattice.

{\it The subcase $a=0$}. Two $D$-invariant lattices such that
$D\bmod (z)$ has the eigenvalues~$\pm \lambda$,
are related by a matrix $\big(\begin{smallmatrix} x_1& x_2\\ x_3& x_4\end{smallmatrix}\big)\in {\rm GL}_2( \mathbb{C}(\!(z)\!))$ satisfying
\[
\left(\begin{matrix} x_1& x_2\\ x_3& x_4\end{matrix}\right) \left(z\frac{{\rm d}}{{\rm d}z}+
\left(\begin{matrix} \lambda& 0\\ 0& -\lambda\end{matrix}\right) \right)
 =\left(z\frac{{\rm d}}{{\rm d}z}+\left(\begin{matrix} \lambda & 0\\ 0& -\lambda\end{matrix}\right) \right)\left(\begin{matrix} x_1& x_2\\ x_3& x_4\end{matrix}\right).
 \]
 The solutions are
$\big(\begin{smallmatrix} x_1& x_2\\ x_3& x_4\end{smallmatrix}\big) =
\Big(\begin{smallmatrix} c_1& c_2z^m\\ c_3z^{-m}& c_4 \end{smallmatrix}\Big)$ with $c_1,c_2,c_3,c_4\in \mathbb{C}$ and
$c_1c_4-c_2c_3\neq 0$. The two $D$-invariant lattices are equal if and only if $c_3=0$. The conclusion is that
there are many $D$-invariant lattices such that $D$ has, modulo $(z)$, eigenvalues $\pm \lambda$. This proves
statement~(1) of the lemma.

 Using the form of the matrix
 $\Big(\begin{smallmatrix} c_1& c_2z^m\\ c_3z^{-m}& c_4\end{smallmatrix}\Big)$, one deduces statement (2) of the lemma. Furthermore, the case $a=0$ corresponds to $N$ is trivial if $\theta$ is even and $N\otimes L$ is trivial if $ \theta$ is odd.

{\it The case $\theta=0$}. Let $\Lambda$ be a $D$-invariant lattice such that $D$ has on $\Lambda/z\Lambda$
only $0$ as eigenvalue. Then $\Lambda$ has a basis $b_1$, $b_2$ such that $D$ becomes the operator
$z\frac{{\rm d}}{{\rm d}z}+\big(\begin{smallmatrix} 0& a\\ 0& 0\end{smallmatrix}\big)$ with $a\in \mathbb{C}$. The module~$N$ is trivial if and only if $a=0$.
As above one verifies that there is only one $D$-invariant lattice with this property.

Finally, if $a=0$ then a basis $b_1$, $b_2$ exists with $D(b_1)=D(b_2)=0$. Hence $D(v)=0$ if and only if
$v\in \mathbb{C}b_1+\mathbb{C}b_2$. \end{proof}

 \begin{Corollary}\label{2.*} Suppose that the parameters satisfy ``${\rm restr}$'' and
 let $\mathcal{V}$ denote the vector bundle $Oe_1+O(-[\infty])e_2$.
 \begin{enumerate}
\item[$(1)$] Consider non exceptional parameters $\{\theta_*\}$, i.e., excluded are
\[ (\theta=0, -\epsilon_0\theta_0+\epsilon_1\theta_1=-1),\qquad (\theta_0=\pm 1, \theta=\pm \theta_1),
\qquad (\theta_1=\pm 1, \theta=\pm \theta_0). \]
For every module $(M,u)\in {\bf S}(\theta_0,\theta_1,\theta_\infty)$ there is a unique connection
 $\nabla\colon \mathcal{V} \rightarrow \Omega([0]+[1]+2[\infty]))\otimes \mathcal{V}$ with local data defined by the parameters and generic fibre $M$. Moreover, the line bundle $Oe_1$ is not invariant under $\nabla$.

 Conversely, suppose that the connection $\nabla\colon \mathcal{V} \rightarrow \Omega([0]+[1]+2[\infty])\otimes \mathcal{V}$ has the
local data defined by the parameters and that the line bundle $Oe_1$ is not invariant under $\nabla$. Then there is a $(M,u)\in {\bf S}(\theta_0,\theta_1,\theta_\infty)$
where $M$ is the generic fibre of $\nabla$.

\item[$(2)$] Consider exceptional parameters $\{\theta_*\}$.
For every $(M,u)\in {\bf S}(\theta_0,\theta_1,\theta_\infty, \epsilon_2=-1)$ there is a unique connection
 $\nabla\colon \mathcal{V} \rightarrow \Omega([0]+[1]+2[\infty])\otimes \mathcal{V}$ with the local data defined by the parameters and with generic
fibre $M$. Moreover, for every invariant saturated line bundle $\mathcal{L}\subset \mathcal{V}$ one has $\mathcal{L}\neq Oe_1$ and $\epsilon_2=-1$.
\end{enumerate}
Conversely, for any connection $\nabla$ with the above properties and generic fibre $M$ there exists an~element
$(M,u)\in {\bf S}(\theta_0,\theta_1,\theta_\infty, \epsilon_2=-1)$.

 The above statements hold with $\epsilon_2=-1$ replaced by $\epsilon_2=1$.
\end{Corollary}

 \begin{proof} {\it Proof of part }(1).
 Write $D$ for the differential on $M$. We provide $M$ with ``invariant lattices''
 $\Lambda$ at the points $0$, $1$, $\infty$. For $z=0$ this means that
 $\Lambda \subset \mathbb{C}(\!(z)\!)\otimes M$ is a free $\mathbb{C}[\![z]\!]$-module of rank 2
 with the properties $\mathbb{C}(\!(z)\!)\otimes \Lambda =\mathbb{C}(\!(z)\!)\otimes M$ and
 $zD(\Lambda)\subset \Lambda$. We require that $zD$ has on $\Lambda/z\Lambda$ eigenvalues
 $\pm \frac{\theta_0}{2}$. The condition on the lattice at $z=1$ is similar. The condition on the lattice
 $\Lambda$ at $z=\infty$ is that $\Lambda$ has a basis such that $D$ on this basis has the
 operator form
 $\frac{{\rm d}}{{\rm d}z}+\frac{1}{z} \Big(\begin{smallmatrix} x-\frac{1}{2}& 0 \\ 0& -x-\frac{1}{2} \end{smallmatrix}\Big)$
 where $x=\frac{tz+\theta}{2}$ and $\theta=\theta_\infty +1$.

 The first two lattices exist and are unique if $\theta_0,\theta_1\not \in \mathbb{Z}$. For
 $\theta_0\in \mathbb{Z}$ and/or $\theta_1\in \mathbb{Z}$ unique lattices are produced by
 Lemma~\ref{2.1}. One easily verifies the existence and uniqueness for the invariant lattice at $z=\infty$.

 This leads to a connection
 $\nabla\colon \mathcal{W}\rightarrow \Omega([0]+[1]+2[\infty])\otimes \mathcal{W}$ with generic
 fiber $M$. Since $\mathcal{W}$ has degree $-1$ we can identify $\mathcal{W}$ with the
 vector bundle $O(k[\infty])e_1+O((-1-k)[\infty ])e_2$.

 For $k\geq 1$, the operator
 $D=\nabla_{\frac{{\rm d}}{{\rm d}z}}$ has on the basis $e_1$, $e_2$ the form
 $\frac{{\rm d}}{{\rm d}z}+\frac{1}{z(z-1)} \big(\begin{smallmatrix} a &c \\ 0&-a \end{smallmatrix}\big) $
 with $a=a_0+a_1z+a_2z^2$ and $c$ a polynomial of some degree $\leq 3+2k$.
 One finds equations
 \[a_0^2=\bigg(\frac{\theta_0}{2}\bigg)^2, \qquad (a_0+a_1+a_2)^2=\bigg(\frac{\theta_1}{2}\bigg)^2 \]
and by using the generators $z^ke_1$, $z^{-1-k}e_2$ of $\mathcal{W}$ at $\infty$ one obtains
 \[ a_2^2=\bigg(\frac{t}{2}\bigg)^2, \qquad 2a_1a_2+2a_2^2+(1+2k)a_2=\frac{t \theta}{2}.\]
 Thus there are $(\epsilon_0,\epsilon_1,\epsilon_2)\in \{\pm 1\}^3$ with
 \begin{gather*}
 a_0=\epsilon_0\frac{\theta_0}{2},\! \qquad a_2=\epsilon_2\frac{t}{2},\! \qquad
 a_1=-\epsilon_0\frac{\theta_0}{2}+\epsilon_1\frac{\theta_1}{2}-\epsilon_2\frac{t}{2},\! \qquad
 k=-\frac{1}{2}+\epsilon_2\frac{\theta}{2}-\epsilon_1\frac{\theta_1}{2}+\epsilon_0\frac{\theta_0}{2}.
 \end{gather*}
 This contradicts part (5) of the assumptions on ${\bf S}(\theta_0,\theta_1,\theta_\infty)$.
 Thus we have $k=0$ and $\mathcal{W}=\mathcal{V}$. If the line bundle $Oe_1$ is invariant
 under $\nabla$, then the same computation provides a contradiction.

 A similar computation proves the last statement of part (1).

{\it Proof of part} (2). It suffices to remark that, by definition, the reducible connections $\nabla$ that are admitted
 correspond with the reducible objects that are admitted in ${\bf S}(\theta_0,\theta_1,\theta_\infty, \epsilon_2=-1)$. The same
 remark holds for $\epsilon_2=1$. \end{proof}

 \begin{Remark}\label{1.3} \rm Let $M$ be a reducible object of ${\bf S}(\theta_0,\theta_1,\theta_\infty)$. The submodule
 $L$ of $M$ is represented by a line bundle $\mathcal{L}\subset \mathcal{V}$ with local data
 $\frac{{\rm d}}{{\rm d}z}+\epsilon_0\frac{\theta_0}{2z}$, $\frac{{\rm d}}{{\rm d}z} -\epsilon_1\frac{\theta_1}{2(z-1)}$ and
 $\frac{{\rm d}}{{\rm d}z}-\epsilon_2\frac{tz+\theta}{2z}-\frac{1}{2z}$. Therefore, $\mathcal{V}/\mathcal{L}$ is again a
 line bundle. By the assumption ``${\rm restr}$'', the degree of $\mathcal{L}$ is $0$ or $-1$. Since degree zero is excluded by part (5) of the definition of ${\bf S}(\theta_0,\theta_1,\theta_\infty)$, the degree of $\mathcal{L}$ is $-1$.
 \end{Remark}

 \section[Definition of the monodromy space R\^{}\{geom\}(theta\_0,theta\_1,theta\_infty)]{Definition of the monodromy space $\boldsymbol{\mathcal{R}^{\rm geom}(\theta_0,\theta_1,\theta_\infty)}$}\label{Section2}
{\it The next step is to convert the analytic data for the modules in ${\bf S}(\theta_0,\theta_1,\theta_\infty)$
into a moduli space. As before, we assume that ``${\rm restr}$'' holds for the parameters. }

 As in Section~\ref{Section1}, additional data are added in order to have a good ``geometric quotient''
without singularities. The {\it analytic data} associated to any $(M,u)\in {\bf S}(\theta_0,\theta_1,\theta_\infty)$ is represented by
the five matrices defined below.
 \begin{enumerate}\itemsep=0pt
\item[$({\rm i})$] The monodromy matrix $\big(\begin{smallmatrix} a_1& b_1\\ c_1& d_1\end{smallmatrix}\big)$
 at $z=0$ with eigenvalues ${\rm e}^{\pm \pi {\rm i} \theta_0}$ and with $s_1={\rm e}^{ \pi {\rm i} \theta_0}+{\rm e}^{- \pi {\rm i} \theta_0}$.
\item[$({\rm ii})$] The monodromy matrix $\big(\begin{smallmatrix} a_2& b_2\\ c_2& d_2\end{smallmatrix}\big)$
at $z=1$ with eigenvalues ${\rm e}^{\pm \pi {\rm i} \theta_1}$ and with $s_2={\rm e}^{\pi {\rm i} \theta_1}+{\rm e}^{- \pi {\rm i} \theta_1}$.
\item[$({\rm iii})$] The formal monodromy matrix
$\big(\begin{smallmatrix} \alpha & 0\\ 0& \frac{1}{\alpha}\end{smallmatrix}\big)$ at $z=\infty$ with eigenvalue $s_3:=\alpha = {\rm e}^{\pi {\rm i} \theta_\infty}$ and the two Stokes matrices
$\big(\begin{smallmatrix} 1& 0\\ f_2& 1\end{smallmatrix}\big) $, $\big(\begin{smallmatrix} 1& f_1\\ 0& 1\end{smallmatrix}\big)$.
\item[$({\rm iv})$] The relation between these matrices reads, with a slight variation on the notation used in \cite[Section~3.2.1]{vdP-Sa},
 \[\left(\begin{matrix} a_1& b_1\\ c_1& d_1\end{matrix}\right) \cdot \left(\begin{matrix}a_2& b_2\\ c_2& d_2\end{matrix}\right) =
\left(\begin{matrix}\alpha& 0\\ 0& \frac{1}{\alpha}\end{matrix}\right)
 \cdot \left(\begin{matrix}1& 0\\ f_2& 1\end{matrix}\right) \cdot \left(\begin{matrix}1& f_1\\ 0& 1\end{matrix}\right).\]

The module $M$ has a formal solution space at $z=\infty$ with a natural basis $\{B_1,B_2\}$ (see below), unique up to multiplication by scalars.
 The five matrices above are taken with respect to this basis and the {\it monodromy space} is, roughly speaking, the solution space of the above system of equations, divided out by the action of $\mathbb{G}_m$ which is
given by conjugation with the matrices $\big\{ \big(\begin{smallmatrix} c & 0\\ 0& 1\end{smallmatrix}\big) \mid c\in \mathbb{C}^*\big\}$.

The split case, i.e., $b_1=b_2=c_1=c_2=f_1=f_2=0$, is $\mathbb{G}_m$-invariant and
produces a~singularity. After leaving out the split case, the remaining problem is that the quotient space could be non-separated.
\label{twee}

The remedy for this, chosen here,
is to construct two geometric quotients. The first one, which we denote by $\mathcal{R}^{\rm geom}(\theta_0,\theta_1,\theta_\infty, \epsilon_2=-1)$,
is obtained by division of the above space under the condition $(c_1,c_2)\neq 0$ (see this section for details).

The second one, denoted by
$\mathcal{R}^{\rm geom}(\theta_0,\theta_1,\theta_\infty, \epsilon_2=1)$, is the geometric quotient under the condition $(b_1,b_2)\neq 0$.


\item[$({\rm v})$] \label{remv} As in part (5) of the definition of ${\bf S}(\theta_0,\theta_1,\theta_\infty)$,
we assume that the parameters satisfy ``${\rm restr}$''. If the parameters
$\theta_0$, $\theta_1$, $\theta_\infty$ allow reducible modules, then there is a solution
${(\epsilon_0,\epsilon_1,\epsilon_2)\in \{\pm 1\}^3}$ of the equation
$ -1=-\frac{1}{2}-\epsilon_2\frac{\theta}{2}+\epsilon_1\frac{\theta_1}{2}-\epsilon_0\frac{\theta_0}{2}$
(see Section~\ref{sec1.1}, Corollary~\ref{2.*} and Remark~\ref{1.3}).
A differential module $M\in {\bf S}(\theta_0,\theta_1,\theta_\infty)$ produces a solution for the matrix equation of~(iv). The
cases where reducible objects are present:
 \begin{enumerate}\itemsep=0pt
\item[$({\rm a})$] Suppose that $\epsilon_2=-1$ holds for all reducible $M\in {\bf S}(\theta_0,\theta_1,\theta_\infty)$.
Then $(c_1,c_2)\neq 0$ holds and the monodromy data define a canonical map
\[{\bf S}(\theta_0,\theta_1,\theta_\infty)\rightarrow \mathcal{R}^{\rm geom}(\theta_0,\theta_1,\theta_\infty, \epsilon_2=-1)(\mathbb{C}).\]

\item[$({\rm b})$] Suppose that $\epsilon_2=1$ holds for all reducible $M\in {\bf S}(\theta_0,\theta_1,\theta_\infty)$.
 Then $(b_1,b_2)\neq 0$ holds and the monodromy data define a canonical map
\[{\bf S}(\theta_0,\theta_1,\theta_\infty)\rightarrow \mathcal{R}^{\rm geom}(\theta_0,\theta_1,\theta_\infty, \epsilon_2=1)(\mathbb{C}).\]
\end{enumerate}

Indeed, suppose that the module $M$ has a unique one-dimensional submodule $L$ and that~$L$ satisfies (v) with
$\epsilon_2=-1$.
The formal module $M\otimes \mathbb{C}(\!(1/z)\!)$ has eigenvalues~$\pm \frac{tz}{2}$ and is classified by
the tuple $\big(V, \big\{V_{\frac{tz}{2}},V_{-\frac{tz}{2}}\big\},\gamma\big)$. Write $B_1$ and $B_2$ for bases of~$V_{\frac{tz}{2}}$ and
$V_{-\frac{tz}{2}}$. Then $\{B_1,B_2\}$ is a basis used for the matrices defining the spaces
$\mathcal{R}^{\rm geom}(\theta_0,\theta_1,\theta_\infty, \epsilon_2=\pm 1)$.
Locally at $z=\infty$, $L$ has the form $\frac{{\rm d}}{{\rm d}z}+\frac{tz+\theta}{2z}-\frac{1}{2z}$ and has (generalized) eigenvalue~$-\frac{tz}{2}$.
This implies that the line $\mathbb{C}B_2\subset V$ is invariant under the monodromy data
and so ${b_1=b_2=f_1=0}$. Since $M$ has no other one-dimensional submodules, one finds $(c_1,c_2)\neq (0,0)$.
This proves statement (a). Case (b) is similar.
 \begin{enumerate}\itemsep=0pt
\item[$({\rm c})$] There are a few cases, called {\it the exceptional parameters} $\{\theta_*\}$, such that both
 $\epsilon_2=-1$ and $\epsilon_2=1$ occurs for reducible modules.
 The first case is $\theta =0$ and $-\epsilon_0\theta_0+\epsilon_1\theta_1=-1$ for some $\epsilon_0,\epsilon_1\in \{\pm 1\}$.
The tuples $(\epsilon_0,\epsilon_1,1)$ and $(\epsilon_0,\epsilon_1,-1)$ describe the same family of modules but with the variable $t$
replaced by $-t$. There is a canonical map ${\bf S}(\theta_0,\theta_1,\theta_\infty, \epsilon_2=-1)\rightarrow \mathcal{R}^{\rm geom}(\theta_0,\theta_1,\theta_\infty, \epsilon_2=- 1)(\mathbb{C})$. The same holds with $\epsilon_2=-1$ replaced by $\epsilon_2=1$.
The other cases are $(\theta_1=\pm 1, \theta=\pm \theta_0)$ and $(\theta_0=\pm 1, \theta=\pm \theta_1)$.
Here however, the two values for $\epsilon_2$ define distinct families of differential modules.
Again, there are canonical maps
\[{\bf S}(\theta_0,\theta_1,\theta_\infty, \epsilon_2=\pm 1)\rightarrow \mathcal{R}^{\rm geom}(\theta_0,\theta_1,\theta_\infty, \epsilon_2=\pm1)(\mathbb{C}).\]

\item[$({\rm d})$] The construction of $\mathcal{R}^{\rm geom}(\theta_0,\theta_1,\theta_\infty, \epsilon_2=-1)$, i.e.,
the monodromy space under the condition of (a), namely $(c_1,c_2)\neq (0,0)$, will be worked out here in some detail. The second case where the condition is $(b_1,b_2)\neq (0,0)$, is similar and will not be worked out.
 {\it In the sequel we will write $\mathcal{R}^{\rm geom}(\theta_0,\theta_1,\theta_\infty)$ if from the context the choice of
$\epsilon_2\in \{\pm 1\}$ is clear or has no importance. Only in the case of exceptional parameters we will keep the above notation.}

\item[$({\rm e})$] The space $\mathcal{R}^{\rm geom}(\theta_0,\theta_1,\theta_\infty)$ depends, apart
from the choice of reducible objects, only on the variables $s_1,s_2,s_3=\alpha $.

\item[$({\rm f})$] According to results in \cite[Section~3.2.2]{vdP-Sa}, the space
$\mathcal{R}^{\rm geom}(\theta_0,\theta_1,\theta_\infty)$ (say, under the condition $(c_1,c_2)\neq 0$) has
 a nonempty reducible locus if and only if $\big(s_3+s_3^{-1}\big)^2-s_1s_2\big(s_3+s_3^{-1}\big)+s_1^2+s_2^2=4$ holds.
For $s_3=\pm 1$, $s_1\neq \pm 2$, $s_2=\pm s_1$, the reducible locus consists of two, non-intersecting, projective lines.
For the other cases with $s_1, s_2\neq \pm 2$ the reducible locus is one projective line.
\end{enumerate}

\item[$({\rm vi})$] If $\theta_0\in \mathbb{Z}$ or equivalently $s_1=\pm 2$, then there is the additional datum of an eigenline $\mathbb{C}v\neq 0$ for the
monodromy matrix $\big(\begin{smallmatrix} a_1& b_1\\ c_1& d_1\end{smallmatrix}\big)$.
\item[$({\rm vii})$] If $\theta_1\in \mathbb{Z}$ or equivalently $s_2=\pm 2$, then there is the additional datum of an eigenline $\mathbb{C}v\neq 0$ for the
monodromy matrix $\big(\begin{smallmatrix} a_2& b_2\\ c_2& d_2\end{smallmatrix}\big)$.
\end{enumerate}

For the explanation of (i)--(vii), we recall from \cite[Section~3.2]{vdP-Sa} the following
description of the map ${\bf S}(\theta_0,\theta_1,\theta_\infty)\rightarrow \mathcal{R}^{\rm geom}(\theta_0,\theta_1,\theta_\infty)$. Let an object $(M,u)$ be given, or
a connection in the (to be constructed) space $\mathcal{M}(\theta_0,\theta_1,\theta_\infty)$.
 At $z=\infty$ the differential module $\mathbb{C}\big(\!\big(z^{-1}\big)\!\big)\otimes M$ is equivalent to the differential operator
$\frac{{\rm d}}{{\rm d}z}-\frac{1}{z}\cdot\Big(\begin{smallmatrix} \frac{tz+\theta_\infty}{2}& 0\\ 0& -\frac{tz+\theta_\infty}{2}\end{smallmatrix}\Big) $. Write
$y=z^{\theta_\infty /2}{\rm e}^{tz/2}$. The solutions $\big\{\big(\begin{smallmatrix}y\\ 0 \end{smallmatrix}\big) , \big(\begin{smallmatrix} 0\\ y^{-1} \end{smallmatrix}\big) \big\}$ of this differential operator produce
a basis $\{B_1,B_2\}$ for the space $V$ of the formal solutions of $M$ at $z=\infty$. In fact, as before, $V_{\frac{tz}{2}}=\mathbb{C}B_1$ and
$V_{\frac{-tz}{2}}=\mathbb{C}B_2$.
The $B_1$ and $B_2$ are unique up to multiplication by
elements of $\mathbb{C}^*$. Furthermore one chooses paths from~$0$ to~$1$ and from~$1$ to~$\infty$. The latter should have a nonsingular direction at $\infty$.

Multisummation at this nonsingular direction lifts $B_1$, $B_2$ to solutions of the differential module $M$ on a large sector around this direction. These extend to multivalued functions
on $\mathbb{P}^1\setminus \{0,1,\infty \}$ and the monodromies at $z=0,1,\infty$
and the Stokes matrices are well defined and are all considered with respect to this basis $B_1$, $B_2$.
For suitable orientations, the product of the topological monodromies at $z=0$ and $z=1$
is equal to the topological monodromy at $z=\infty$. Then the identity (iv) follows from the ``monodromy identity'' \cite[Proposition 8.12]{vdP-S}. The above five matrices are considered modulo simultaneous conjugation by $\mathbb{G}_m=\big\{\big(\begin{smallmatrix} c& 0\\ 0& 1 \end{smallmatrix}\big) \mid c\in \mathbb{C}^*\big\}$,
 since the basis $B_1$, $B_2$ is unique up to multiplication by scalars.

 \begin{Remark}\label{2.10} For a single module $M$, the definition of the image of $M$
 into $\mathcal{R}^{\rm geom}(\theta_0,\theta_1,\theta_\infty)$ poses no problems since one can
 always choose a nonsingular direction for the multisummation.

 However for a family this
 poses a problem since the two singular directions
 $-\beta, \frac{1}{2}-\beta \in \mathbb{R}/\mathbb{Z}$ vary with
 $\frac{t}{|t|}={\rm e}^{2\pi {\rm i} \beta}$. The remedy is to choose a nonsingular direction,
 say $-\beta +\frac{1}{4}$, that moves with $\frac{t}{|t|}$. This can be done, in a continuous way, by using $u$ with $t={\rm e}^{2\pi {\rm i} u}$ on the universal covering of $\mathbb{P}^1\setminus \{0,1,\infty \}$.
\end{Remark}

{\it Construction of $\mathcal{R}^{\rm geom}(\theta_0,\theta_1,\theta_\infty)$, separation and parabolic structure}.
Let $\operatorname{Spec}(R)$ denote the affine variety describing the above five matrices and their relations (for fixed
$s_1:={\rm e}^{ \pi {\rm i} \theta_0} +{\rm e}^{- \pi {\rm i} \theta_0}$, $s_2:={\rm e}^{ \pi {\rm i} \theta_1} +{\rm e}^{-\pi {\rm i} \theta_1}$, $s_3=\alpha ={\rm e}^{\pi {\rm i} \theta_\infty}$). {\it We forget for the moment restrictions related to reducibility}.
After elimination of $f_1$ and $f_2$ one finds that $R=\mathbb{C}[a_1,b_1,c_1,d_1,a_2,b_2,c_2,d_2]/\mbox{(relations)}$.
The relations are generated by
\begin{gather*}
a_1d_1-b_1c_1=1, \quad a_2d_2-b_2c_2=1, \quad a_1+d_1=s_1, \quad a_2+d_2=s_2, \quad a_1a_2+b_1c_2=s_3.
\end{gather*}
One needs a quotient of $\operatorname{Spec}(R)$ by the action of $\mathbb{G}_m$.
 A categorical quotient $\mathcal{R}(\theta_0,\theta_1,\theta_\infty):=\operatorname{Spec}(R_0)$ with
$R_0=R^{\mathbb{G}_m}$ is computed in \cite[Section~3.2]{ vdP-Sa}. It is an affine cubic surface with three lines at
 infinity. This is given an interpretation in \cite{C-M-R}. The disadvantage is that it is not a geometric
quotient and, moreover, there are many singularities (depending on $s_1$, $s_2$ and~$s_3$). These singularities are caused
by resonance, reducibility and fixed points for the action of $\mathbb{G}_m$ (see the table in \cite[Section~3.2.2]{vdP-Sa}). Here we will define a smooth ``geometric quotient''
$\mathcal{R}^{\rm geom}(\theta_0,\theta_1,\theta_\infty)$ by assuming
 (v) (i.e., we assume $(c_1,c_2)\neq (0,0)$) and adding additional data for $s_1=\pm 2$ and/or $s_2=\pm 2$ (see~(vi) and~(vii)).

 {\it Suppose that $\theta _0,\theta _1\not \in \mathbb{Z}$ or equivalently $s_1\neq \pm 2$, $s_2\neq \pm 2$}.
 Then we consider the open subspace $X$ of $\operatorname{Spec}(R)$, defined by condition (v), namely
 the case $(c_1,c_2)\neq (0,0)$. A geometric quotient $X/\mathbb{G}_m$ exists and is connected, smooth and separated as can be seen as follows.

 The space $X$ is the union of the two Zariski open affine subspaces given by $c_1\neq 0$ and~$c_2\neq 0$. Let $X(c_1\neq 0)\subset X\subset \operatorname{Spec}(R)$ denote the affine subspace given by $c_1\neq 0$ and let $X({c_1\! =\! 1})$ be the closed subset given by $c_1=1$. The action of $\mathbb{G}_m$ on $X$ yields an isomorphism
$\mathbb{G}_m\times X({c_1\!=\!1})\rightarrow X(c_1\neq 0)$. Therefore, we can identify
$X(c_1\neq 0)/\mathbb{G}_m$ with $X(c_1\!=\!1)$.
The affine space $X(c_1\!=\!1)$ is given by variables $a_1$, $b_1$, $d_1$, $a_2$, $b_2$, $c_2$, $d_2$ and relations
\[
d_j=s_j-a_j, \ j=1,2,\qquad a_1d_1-b_1=1,\qquad a_2d_2-b_2c_2=1, \qquad a_1a_2+b_1c_2=s_3\neq 0.
\]
After elimination $d_1$, $d_2$, $b_1$ the coordinate ring of $X(c_1=1)$ is equal to
$\mathbb{C}[a_1,a_2,b_2,c_2]/(a_1a_2+(-1+a_1(s_1-a_1)c_2-s_3, a_2(s_1-a_2)-b_2c_2-1))$. A calculation shows that
$X(c_1=1)$ is connected, has dimension 2 and is smooth under the assumption that $s_2\neq \pm 2$.

A similar calculation shows that $X(c_2\neq 0 )\cong X(c_2=1)\times \mathbb{G}_m$ and that $X(c_2=1)$ has dimension 2,
is connected and smooth under the assumption that $s_1\neq \pm 2$.

Let $X/\mathbb{G}_m$ denote the gluing of $X(c_1\neq 0)/ \mathbb{G}_m$ and $X(c_2\neq 0)/\mathbb{G}_m$. In order to
verify that $X/\mathbb{G}_m$ is {\it separated} we consider the morphism $F\colon X\rightarrow \mathbb{P}^1$ that sends a point of $X$ to the
equivalence class $[c_1:c_2]\in \mathbb{P}^1$ of $(c_1,c_2)$. Since $F$ is $\mathbb{G}_m$-equivariant and $\mathbb{P}^1$ has trivial
$\mathbb{G}_m$-action, one obtains a morphism $G\colon X/\mathbb{G}_m\rightarrow \mathbb{P}^1$. The fibers of $F$ are two-dimensional affine
spaces~$V$ with a $\mathbb{G}_m$-action such that $V/\mathbb{G}_m$ is affine and $V/\mathbb{G}_m \times \mathbb{G}_m\cong V$. Then the fibers of~$G$ are the affine one-dimensional spaces $V/ \mathbb{G}_m$. We apply the valuative criterion of separatedness
(see \cite[Theorem II.4.3]{Hart}) to $X/\mathbb{G}_m$. Since $\mathbb{P}^1$ is separated it suffices to verify that the fibers of~$G$
are separated. The latter holds since these fibers are affine.

 For $s_1\neq \pm 2$, $s_2\neq \pm 2$, the space $\mathcal{R}^{\rm geom}(\theta_0,\theta_1,\theta_\infty)$ is by definition the above space $X/\mathbb{G}_m$.

{\it The cases $\theta_0\in \mathbb{Z}$ or $\theta_1\in \mathbb{Z}$. We explain the case
 $\theta_0=2m\in 2\mathbb{Z}$ in some detail}. The other cases $\theta_0\in 1+2\mathbb{Z}$ or $\theta_1\in \mathbb{Z}$
 can be treated in a similar way. The space $\mathcal{R}^{\rm geom}(2m,\theta_1,\theta_\infty)$ is given by the tuples (modulo $\mathbb{G}_m$-action)
\[ \left(\begin{pmatrix} a_1& b_1\\ c_1& d_1\end{pmatrix}, \begin{pmatrix} a_2& b_2\\ c_2& d_2\end{pmatrix}, [y_1:y_2] \right)\in \mathbb{A}^4\times \mathbb{A}^4\times \mathbb{P}^1\]
with the following relations:
\begin{gather*}
 a_1+d_1=2 , \qquad a_1d_1-b_1c_1=1 , \qquad a_2+d_2=s_2 ,\qquad a_2d_2-b_2c_2=1 ,\\ a_1a_2+b_1c_2=s_3 ,\qquad
\begin{pmatrix} a_1& b_1\\ c_1& d_1\end{pmatrix} \begin{pmatrix} y_1\\ y_2\end{pmatrix} =\begin{pmatrix} y_1\\ y_2\end{pmatrix},
\end{gather*}
 and one of the elements $c_1,c_2 \neq 0$. We consider again the geometric quotient $\mathcal{R}^{\rm geom}(2m,\theta_1,\theta_\infty)$ of the above space by the action of $\mathbb{G}_m$. This quotient is obtained by gluing some affine open subspaces. Each subspace is obtained by normalizing one of the $c_1$, $c_2$ to $1$ and also one of the elements $y_1$, $y_2$ is normalized to 1.
For each of these four open affine subsets one can make the
coordinate ring explicit and verify that there are no singularities.

 As an example we make the case $c_2=1$, $y_2=1$ explicit. Elimination of the variables $a_1$, $d_1$, $d_2$, $b_1$, $b_2$ leads to the
 nonsingular coordinate ring $\mathbb{C}[a_2,c_1,y_1]/\big(a_2c_1y_1-c_1y_1^2+a_2-s_3\big)$ of this affine subset.

 As in the case $s_1\neq \pm 2$, $s_2\neq \pm 2$ one can show that the gluing of the four affine spaces produces a connected, separated and
 smooth variety.

\begin{Proposition}\label{2.2}\label{prop2.1} Suppose that the parameters satisfy ``${\rm restr}$''.
 Let $s_1={\rm e}^{\pi {\rm i} \theta_0}+{\rm e}^{-\pi {\rm i} \theta_0}$,
$s_2={\rm e}^{\pi {\rm i} \theta_1}+{\rm e}^{-\pi {\rm i} \theta_1}$, $s_3={\rm e}^{\pi {\rm i} \theta_\infty}$.
\begin{enumerate}\itemsep=0pt
\item[$(1)$] $\mathcal{R}^{\rm geom}(\theta_0,\theta_1,\theta_\infty)$ is a simply connected, smooth variety of dimension $2$. It is a geometric quotient
of the analytic data including the eigenlines for the cases $s_1=\pm 2$ and/or $s_2=\pm 2$.
\item[$(2)$] Suppose that the parameters are not exceptional.
The map
\[{\bf S}(\theta_0,\theta_1,\theta_\infty)\rightarrow \mathcal{R}^{\rm geom}(\theta_0,\theta_1,\theta_\infty)(\mathbb{C})\times \mathbb{C},\]
 which associates to any $(M,u)$ the analytic data of $M$ {\rm (}including additional lines for $\theta_0\in \mathbb{Z}$ and/or $\theta_1\in \mathbb{Z}${\rm)} and $u\in \mathbb{C}$, is a bijection.
\item[$(3)$] Suppose that the parameters are exceptional. Then {\rm (2)} remains valid after adding the data $\epsilon_2=-1$ or $\epsilon_2=1$
to both sets.
\end{enumerate}
 \end{Proposition}

\begin{proof} Most of the proof of (1) is sketched above.
For the proof that the complex algebraic variety $\mathcal{R}^{\rm geom}(\theta_0,\theta_1,\theta_\infty)$ is {\it simply connected} one considers
the morphism $\operatorname{pr}\colon \mathcal{R}^{\rm geom}(\theta_0,\theta_1,\theta_\infty)\rightarrow \mathbb{C}^2$ which sends
the data $\big(\begin{smallmatrix} a_1& b_1\\ c_1& d_1 \end{smallmatrix}\big)$,
$\big(\begin{smallmatrix} a_2& b_2\\ c_2& d_2\end{smallmatrix}\big) $ and $[y_1:y_2]$ (for special $s_1$, $s_2$)
to $(a_1,a_2)\in \mathbb{C}^2$.

First we consider the case $s_1\neq \pm 2$ and $s_2\neq \pm 2$.
 A computation shows that the fibre of $pr$ consists of 1 point if $s_3-a_1a_2\neq 0$. For the case
 $s_3-a_1a_2=0$, $a_1(s_1-a_1)-1\neq 0$, $a_2(s_2-a_2)-1\neq 0$ the fibre is empty. For the cases
 $s_3-a_1a_2=0$, $a_1(s_1-a_1)-1\neq 0$, $a_2(s_2-a_2)-1= 0$ and $s_3-a_1a_2=0$, $a_1(s_1-a_1)-1=0$, $a_2(s_2-a_2)-1\neq 0$
 the fibre is $\mathbb{C}$. Finally, for the case $s_3-a_1a_2=0$, $a_1(s_1-a_1)-1= 0$, $a_2(s_2-a_2)-1= 0$, the fibre is the union of three
 complex lines $L_1\cup L_2\cup L_3$. The only intersections are the points $\{p_{1,2}\}=L_1\cap L_2$ and $\{p_{2,3}\}=L_2\cap L_3$
 and $p_{1,2}\neq p_{2,3}$. This shows that the image of $pr$ and that the nonempty fibres of $pr$ are simply connected. Thus by an argument in the spirit of Van Kampen's theorem
 it follows that
 $\mathcal{R}^{\rm geom}(\theta_0,\theta_1,\theta_\infty)$ is simply connected. In the cases $s_1=\pm 2$ and/or $s_2=\pm 2$ the image of $pr$ is the same
 as above. Some of the nonempty fibres have acquired, in comparison with above, projective lines over $\mathbb{C}$ as components.
 We conclude again that $\mathcal{R}^{\rm geom}(\theta_0,\theta_1,\theta_\infty)$ is simply connected.

 The proof of (2) follows by combining the following observations:
 \begin{enumerate}\itemsep=0pt
\item[$({\rm a})$] Given formal and analytic data determine a unique differential module over $\mathbb{C}(z)$,
 see \cite[Theorem 1.7]{vdP-Sa}.
\item[$({\rm b})$] $\mathcal{R}^{\rm geom}(\theta_0,\theta_1,\theta_\infty)$ is by construction a geometric quotient and thus
 $\mathcal{R}^{\rm geom}(\theta_0,\theta_1,\theta_\infty)(\mathbb{C})$ is equal as a set to the analytic data.
\item[$({\rm c})$] ${\bf S}(\theta_0,\theta_1,\theta_\infty)$ and
$\mathcal{R}^{\rm geom}(\theta_0,\theta_1,\theta_\infty)$ are
provided with the same additional data for the cases $\theta_0\in \mathbb{Z}$ and/or $\theta_1\in \mathbb{Z}$.
\end{enumerate}
 The additional remark for the proof of (3) is that the admitted reducible
cases for both sets correspond. \end{proof}

\begin{Remark}\label{remark3.3} We consider the obvious morphism
 \[r\colon \ \mathcal{R}^{\rm geom}(\theta_0,\theta_1,\theta_\infty)\rightarrow
 \mathcal{R}(\theta_0,\theta_1,\theta_\infty)\]
 for the case that the first space is constructed with the condition $(c_1,c_2)\neq (0,0)$. As before
 $\mathcal{R}(\theta_0,\theta_1,\theta_\infty)$ denotes the categorical quotient considered in \cite[Sections~3.2 and~4.3]{vdP-Sa}.
\begin{enumerate}\itemsep=0pt
\item[$({\rm a})$] For general parameters, i.e., no resonance and no reducible data,
the map $r$ is an isomorphism.

\item[$({\rm b})$] 
Suppose the $\theta_*$ satisfy ``${\rm restr}$'' and are not exceptional. Then $r$ is a minimal resolution, as can be seen as follows.
A reducible locus in $\mathcal{R}^{\rm geom}(\theta_0,\theta_1,\theta_{\infty})$ produces a $\mathbb{P}^1$ in the fibre of $r$.
The same holds for the cases $\theta_0\in\mathbb{Z}$ and/or $\theta_1\in\mathbb{Z}$. The assumption ``not exceptional'' implies
that the fibres are proper. By counting the number of $\mathbb{P}^1$'s in the fibres and comparing with the table in
\cite[Table~3.1]{vdP-Sa}, we conclude that the resolution $r$ is minimal.

\item[$({\rm c})$]
 In all cases the map $r$ is birational. For some examples of exceptional parameters, the fibres of $r$ are not complete.
\end{enumerate}
\end{Remark}
\section[Construction of M(theta\_0,theta\_1,theta\_infty)]{Construction of $\boldsymbol{\mathcal{M}(\theta_0,\theta_1,\theta_\infty)}$} \label{Section3}
{\it The next step is a subtle construction of moduli spaces of connections}.
First, we give a construction of this moduli space of connections without paying attention to a parabolic structure of eigenlines. For $\theta_0, \theta_1\not \in \mathbb{Z}$ , this is sufficient. For $\theta_0\in \mathbb{Z}$
and/or $\theta_1\in \mathbb{Z}$, we will refine this construction since the space $\mathcal{M}(\theta_0,\theta_1,\theta_\infty)$ can otherwise have singularities and $\mathcal{M}(\theta_0,\theta_1,\theta_\infty)(\mathbb{C})$ may not coincide with ${\bf S}(\theta_0,\theta_1,\theta_\infty)$.

 {\it The data for the moduli space of connections $\mathcal{M}(\theta_0,\theta_1,\theta_\infty)$ are}
 \begin{enumerate}\itemsep=0pt
\item[$({\rm i})$] The sub-bundle $\mathcal{V}=Oe_1\oplus O(-[\infty])e_2$ of the free bundle $Oe_1\oplus Oe_2$.
\item[$({\rm ii})$]Local connections $L_0$, $L_1$, $L_\infty$, written as differential operators, namely
\begin{enumerate}\itemsep=0pt
\item[$({\rm a})$] $L_0=\frac{{\rm d}}{{\rm d}z}+\frac{1}{z}A_0$ with $A_0$ a $(2\times 2)$-matrix with entries in $\mathbb{C}[\![z]\!]$ and such that
$A_0 \!\! \mod (z)$ has eigenvalues $\pm \frac{\theta_0}{2}$.
\item[$({\rm b})$] $L_1=\frac{{\rm d}}{{\rm d}z}+\frac{1}{z-1}A_1$ with $A_1$ a $(2\times 2)$-matrix with entries in $\mathbb{C}[\![z-1]\!]$ and such that
$A_1\!\! \mod (z-1)$ has eigenvalues $\pm\frac{\theta_1}{2}$.
\item[$({\rm c})$] $L_\infty = \frac{{\rm d}}{{\rm d}z}+
\frac{1}{z}\big(\begin{smallmatrix} x-1/2& 0\\ 0& -x-1/2\end{smallmatrix}\big)$ at $z=\infty$ with $x=\frac{tz+\theta}{2}$,
 $\theta =\theta_\infty +1$, $u\in \mathbb{C}$ and $t={\rm e}^{2\pi {\rm i} u}$.
\end{enumerate}
\end{enumerate}
The {\it moduli functor} $\mathcal{MF}$ (from $\mathbb{C}$-algebras to sets) associates to each $\mathbb{C}$-algebra $R$ the set of equivalence classes of pairs $(\nabla ,u)$ consisting of an element $u\in \mathbb{C} $ and a
 connection
\[\nabla \colon \ R\otimes \mathcal{V}\rightarrow \Omega_{\mathbb{P}^1}([0]+[1]+2\cdot [\infty])\otimes (R\otimes \mathcal{V})\]
satisfying the conditions stated below.
We observe that $D=\nabla_{\frac{{\rm d}}{{\rm d}z}}$ has, with respect to the basis $e_1$, $e_2$, the matrix
 \[ \frac{1}{z(z-1)}\cdot \left(\begin{matrix} a_0+a_1z+a_2z^2 & c_0+c_1z+\cdots +c_3z^3\\ b_0+b_1z& d_0+d_1z+d_2z^2\end{matrix}\right) \]
 with $a_*,b_*,c_*,d_*\in R$.
 The conditions are: $Rb_0+Rb_1=R$ and the action of $D$ on the completions of $R\otimes \mathcal{V}$
 at the points $z=0,1,\infty$ is isomorphic to $R\otimes L_0$, $R\otimes L_1$, $R\otimes L_\infty$, respectively.

Two pairs $(\nabla_1,u_1)$, $(\nabla_2,u_2)$ are called {\it equivalent} if $u_1=u_2$ and there is an
$R$-linear automorphism $g$ of $R\otimes \mathcal{V}$ which transforms $\nabla_1$ into $\nabla_2$.

The $R$-linear automorphisms $g$ of $R\otimes \mathcal{V}$ are given by the formulas
 \[ g(e_1)=\lambda_1e_1, \qquad g(e_2)=\lambda_2e_2+ (\alpha+\beta z)e_1 \qquad\mbox{with}\
\lambda_1, \lambda _2\in R^* \ \mbox{and} \ \alpha, \beta \in R.\]

Our aim is to show that the above functor $\mathcal{MF}$ is representable. The resulting fine moduli space will be denoted by $\mathcal{M}(\theta_0,\theta_1,\theta_\infty)$. Moreover, the map $\mathcal{M}(\theta_0,\theta_1,\theta_\infty)(\mathbb{C})\rightarrow {\bf S}(\theta_0,\theta_1,\theta_\infty)$, which associates to any pair $(\nabla, u)$ the pair $(M,u)$ with
$M$ the generic fibre of $\nabla$, will be shown to be bijective under the assumptions
 $\theta_0,\theta_1\not \in \mathbb{Z}$.

The \label{DefM} above data produce an affine space $M(\theta_0,\theta_1,\theta_\infty)$ given by the variables $t$, $a_0$, $a_1$, $a_2$, $b_0$, $b_1$, $c_0$, $c_1$, $c_2$, $c_3$, $d_0$, $d_1$, $d_2$
and the relations:
 $d_0=-a_0$, $d_1=-a_1$, $d_2=-a_2$ derived from the traces of $L_0$, $L_1$, $L_\infty$, and
 from the determinants of $L_0$, $L_1$, $L_\infty$ one obtains the relations
\begin{alignat*}{3}
& a_0^2+b_0c_0=\frac{\theta_0^2}{4}, \qquad && (a_0+a_1+a_2)^2+(b_0+b_1)(c_0+c_1+c_2+c_3)=\frac{\theta_1^2}{4},&\\
& a_2^2+b_1c_3=\frac{t^2}{4}, \qquad && 2a_2^2+2a_1a_2+a_2+b_0c_3+b_1c_2+2b_1c_3=\frac{t\theta}{2}.&
\end{alignat*}

We note that in case $\theta :=\theta_\infty +1\neq 0$, the variable $t$ can be eliminated by one of the relations.
For $\theta=0$, we will need the variable $t$ for the construction of a good moduli space $\mathcal{M}(\theta_0,\theta_1,\theta_\infty)$.

From the open subspaces of $M(\theta_0,\theta_1,\theta_\infty)$, defined by $b_0\neq 0$ or by $b_1\neq 0$, one has to take the geometric quotient by the action of the automorphism group $G$ of the vector bundle $\mathcal{V}$. Let $M_1(\theta_0,\theta_1,\theta_\infty)$ denote the open subset given
by $b_1\neq 0$. Consider the closed subset $\mathcal{M}_1(\theta_0,\theta_1,\theta_\infty)$ of
$M_1(\theta_0,\theta_1,\theta_\infty)$ defined by $b_1=1$, $a_1=a_2=0$. The action of $G$ defines a
morphism
\[ G\times \mathcal{M}_1(\theta_0,\theta_1,\theta_\infty) \rightarrow M_1(\theta_0,\theta_1,\theta_\infty).\]
This is an isomorphism and therefore we identify $M_1(\theta_0,\theta_1,\theta_\infty)/G$ with
$\mathcal{M}_1(\theta_0,\theta_1,\theta_\infty)$.
 In a~similar way, let $M_2(\theta_0,\theta_1,\theta_\infty)$
be the open subset of $M(\theta_0,\theta_1,\theta_\infty)$ defined by $b_0\neq 0$. Therefore, $M_2(\theta_0,\theta_1,\theta_\infty)/G$ can be identified with $\mathcal{M}_2(\theta_0,\theta_1,\theta_\infty)$
defined as the closed subset of $M(\theta_0,\theta_1,\theta_\infty)$, given by $b_0=1$, $a_0=a_1=0$.
The above construction leads to
\begin{Proposition}\label{2.3} The space $\mathcal{M}(\theta_0,\theta_1,\theta_\infty)$ obtained by gluing
the spaces $\mathcal{M}_1(\theta_0,\theta_1,\theta_\infty)$ and $\mathcal{M}_2(\theta_0,\theta_1,\theta_\infty)$ represents
the functor $\mathcal{MF}$.
\end{Proposition}

{\it More details for the affine space $\mathcal{M}_1=\mathcal{M}_1(\theta_0,\theta_1,\theta_\infty)$}.
 \label{Dmatrix}
Using an automorphism of $\mathcal{V}$, the matrix of $D=\nabla _{\frac{{\rm d}}{{\rm d}z}}$ with respect to basis $e_1$, $e_2$ can uniquely be written as
\begin{gather*}
\frac{{\rm d}}{{\rm d}z}+\frac{1}{z(z-1)}\cdot \left(\begin{matrix} a_0 &c_0+c_1z+c_2z^2+c_3z^3 \\ z+b_0 &-a_0 \end{matrix}\right),
\end{gather*}
with equations
\begin{gather*}
a_0^2+b_0c_0=\frac{\theta_0^2}{4},\qquad a_0^2+(1+b_0)(c_0+c_1+c_2+c_3)=\frac{\theta_1^2}{4},\\
c_3=\frac{t^2}{4},\qquad (b_0+2)c_3+c_2=\frac{\theta t}{2} .
\end{gather*}
Consecutively eliminating $c_3=\frac{t^2}{4}$ and then $c_2=\frac{\theta t}{2}-(b_0+2)\frac{t^2}{4}$, and finally
$c_0=\frac{\theta_1^2}{4}-\frac{\theta_0^2}{4}-{(b_0+1)}\big(c_1-(b_0+1)\frac{t^2}{4}+\frac{\theta t}{2}\big)$
 results in the equation
\[ a_0^2= -\frac{t^2}{4}b_0^3+\biggl(-\frac{t^2}{2}+\frac{\theta t}{2}+c_1\biggr)b_0^2+
\biggl(-\frac{\theta_1^2}{4}+\frac{\theta_0^2}{4}+\frac{\theta t}{2}-\frac{t^2}{4}+c_1\biggr)b_0+\frac{\theta_0^2}{4} \]
for the affine space $\mathcal{M}_1$.
The singular locus of this three-dimensional affine variety over $\mathbb{C}$ is the following:
\begin{enumerate}\itemsep=0pt
\item[$({\rm a})$] For $\theta_0=0$, the line $a_0=b_0=0$ and $c_1=\frac{\theta_1^2}{4}-\frac{\theta t}{2}+\frac{t^2}{4}$.
After substitution of these values in $D$ the singularity at $z=0$ disappears because $c_0=0$.
On the other hand, if the module~$M$ corresponding to a point of $\mathcal{M}_1$ is locally trivial at $z=0$
then one concludes $\theta_0=0$ and $a_0=b_0=c_0=0$.
The singularity at the surface $\mathcal{M}_1$ for $\theta_0=0$ can be resolved by adding as extra data for the
family an invariant line locally at $z=0$. Then the singular point is replaced by a $\mathbb{P}^1$ of directions for the local invariant line.

\item[$({\rm b})$] For $\theta_1=0$, the line $a_0=0$, $b_0=-1$ and $c_1=\frac{\theta_0^2}{4}-\frac{\theta t}{2}$.
It has the same interpretation as in (a). These singular loci coincide with the loci where the map
$\mathcal{M}_1(\theta_0,\theta_1,\theta_\infty)\rightarrow \operatorname{Spec}\big(\mathbb{C}\big[t,t^{-1}\big]\big)$ is not smooth.
\end{enumerate}

{\it More details for $\mathcal{M}_2=\mathcal{M}_2(\theta_0,\theta_1,\theta_\infty)$}.
This Zariski open subset of $\mathcal{M}(\theta_0,\theta_1,\theta_\infty)$ has (apart from $t$) the form
 \[ D=\frac{{\rm d}}{{\rm d}z}+\frac{1}{z(z-1)}\cdot \left(\begin{matrix} A_2z^2 & C_0+\dots +C_3z^3 \\ 1+B_1z & -A_2z^2\end{matrix}\right) ,\]
where
\begin{gather*}
C_0=\frac{\theta_0^2}{4},\qquad A_2^2+(1+B_1)(C_0+C_1+C_2+C_3)=\frac{\theta_1^2}{4},\\
A_2^2+B_1C_3=\frac{t^2}{4},\qquad 2A_2^2+A_2+B_1C_2+2B_1C_3+C_3=\frac{\theta t}{2}, \qquad \theta=\theta_\infty +1.
\end{gather*}
Now $\mathcal{M}\setminus \mathcal{M}_1$ is the closed subset $B_1=0$ of $\mathcal{M}_2$ and has equations
 \[ A_2^2=\frac{t^2}{4},\qquad C_0=\frac{\theta_0^2}{4},\qquad 2A_2^2+A_2+C_3=\frac{\theta t}{2},
 \qquad A_2^2+C_0+C_1+C_2+C_3=\frac{\theta_1^2}{4}.\]
 This defines the disjoint union of two affine lines over $\operatorname{Spec}\big(\mathbb{C}\big[t,t^{-1}\big]\big)$.

\begin{Theorem}\label{2.4} Assume that the parameters satisfy ``${\rm restr}$'' and that
 $\theta_0, \theta_1\not \in \mathbb{Z}$. If the parameters are not exceptional, then the extended Riemann--Hilbert map
${\rm RH}^{+}\colon\mathcal{M}(\theta_0,\theta_1,\theta_\infty)\rightarrow
\mathcal{R}^{\rm geom}(\theta_0,\theta_1,\theta_\infty)\times \mathbb{C}$
is an analytic isomorphism.
If the parameters are exceptional, then the above remains valid after adding $\epsilon_2=-1$ or $\epsilon_2=1$
to both spaces.
 \end{Theorem}
 \begin{proof}
 The morphism of the theorem is a locally defined map in terms of the variable
$t\in \mathbb{C}^*$. By replacing $t$ by $u\in \mathbb{C}$ with $t={\rm e}^{2\pi {\rm i} u}$, it is globally
 defined (see the Remark~\ref{2.10}).
 It is well known that the monodromy matrices and the Stokes matrices depend in an analytic way on the data
 of the connection \cite[Proposition 12.20]{vdP-S}). Thus ${\rm RH}^{+}$ is an analytic map.

 Both spaces are smooth complex algebraic varieties,
 since $\theta_0,\theta_1 \not \in \mathbb{Z}$. It is known that a~bijective analytic map between smooth complex varieties is an analytic isomorphism. Thus we have only to show that ${\rm RH}^{+}$ is bijective.
The bijectivity follows from Proposition~\ref{prop2.1} and Corollary~\ref{2.*}.
\end{proof}

 Theorem~\ref{2.4} implies that the solutions of the fifth Painlev\'e equation (at least for $\theta_0,\theta_1\not \in \mathbb{Z}$) are multivalued meromorphic functions in $t\in \mathbb{C}^*$ or, equivalently,
 meromorphic functions of $u\in \mathbb{C}$ (compare \cite{vdP-T}). In other words, ${\rm P}_5$ has the Painlev\'e property. This was already known from the paper \cite{Mi}.

 It is interesting to note that $\mathcal{M}(\theta_0,\theta_1,\theta_\infty)$ can be identified with the Okamoto--Painlev\'e variety for the fifth Painlev\'e equation (at least for $\theta_0,\theta_1\not \in \mathbb{Z}$
and in the absence of the reducible locus). For the verification, the result of Proposition~\ref{prop2.1} stating that $\mathcal{R}^{\rm geom}(\theta_0,\theta_1,\theta_\infty)$ is {\it simply connected}, is needed.

 \section[The reducible locus in M(theta\_0,theta\_1,theta\_infty)]{The reducible locus in $\boldsymbol{\mathcal{M}(\theta_0,\theta_1,\theta_\infty)}$}\label{Section4}
 We recall that the parameters are supposed to satisfy ``${\rm restr}$''.
 For reducible $(\nabla , u)$, there is a line bundle $\mathcal{W}\subset \mathcal{V}=Oe_1\oplus O(-[\infty ])e_2$, ``invariant'' under $\nabla$ and such that $\mathcal{V}/\mathcal{W}$ is again a line bundle. By construction $\mathcal{W}\neq Oe_1$ and by assumption ``${\rm restr}$'' one can choose~$e_2$ such that $\mathcal{W}=O(-[\infty])e_2$. Furthermore, this $e_2$ is unique up to multiplication by elements in~$\mathbb{C}^*$. Therefore, the set of the closed points of $\mathcal{M}(\theta_0,\theta_1,\theta_\infty)$, corresponding to reducible~$(\nabla ,u)$, consists of the differential operators having with respect to to the basis $\{e_1,e_2\}$ the form
 \[\frac{{\rm d}}{{\rm d}z}+\frac{1}{z(z-1)}\cdot \left(\begin{matrix} a_0+a_1z+a_2z^2 & 0\\ b_0+b_1z & -a_0-a_1z-a_2z^2\end{matrix}\right) \]
with $b_0+b_1z\neq 0$,
considered modulo the action of $\mathbb{G}_m$. The equations are
\[a_0^2=\frac{\theta_0^2}{4},\qquad (a_0+a_1+a_2)^2=\frac{\theta_1^2}{4},\qquad a_2^2=\frac{t^2}{4},\qquad 2a_2^2+2a_1a_2+a_2=\frac{t\theta}{2}.\]
The group $\mathbb{G}_m$ acts trivially on $a_0$, $a_1$, $a_2$. One has
\[a_0=\epsilon_0\frac{\theta_0}{2},\qquad a_1=\epsilon_1\frac{\theta_1}{2}-\epsilon_0\frac{\theta_0}{2}-\epsilon_2\frac{t}{2},\qquad a_2=\epsilon_2 \frac{t}{2} \]
with $\epsilon_0,\epsilon_1,\epsilon_2\in \{1,-1\}$ and $-1=-\frac{1}{2}-\epsilon_2\frac{\theta}{2}+\epsilon_1\frac{\theta_1}{2}-\epsilon_0\frac{\theta_0}{2}$
 (see Remark~\ref{1.3}). Note that there are no conditions on $(b_0,b_1)\neq (0,0)$.

The reducible locus consists of a number of projective lines over $\mathbb{C}\big[t,t^{-1}\big]$ in
$\mathcal{M}(\theta_0,\theta_1,\theta_\infty)$. In general, e.g., if $1$, $\theta_0$, $\theta_1$, $\theta$ are linearly independent over $\mathbb{Q}$, the reducible locus is empty.

We fix a choice of the $\epsilon_j$, $j=0,1,2$. Isomonodromy is obtained by considering $b_0$, $b_1$ as functions of $t$ and completing the Lax pair (compare Section~\ref{Section6}) with a suitable differential operator $\frac{{\rm d}}{{\rm d}t}+ \big(\begin{smallmatrix}\alpha & 0 \\ \beta & -\alpha\end{smallmatrix}\big)$. This results in a Riccati equation for the function $\frac{b_1}{b_0}$. All Riccati solutions for ${\rm P}_5$ are obtained in this way, since one can verify that we found precisely the data of \cite[Theorem~2.1\,(8)]{OO2}. We note that this includes cases with $\theta_0\in \mathbb{Z}$ and/or $\theta_1\in \mathbb{Z}$.

\subsection[Removing split differential modules from the space M(theta\_0,theta\_1,theta\_infty)]{Removing split differential modules from the space $\boldsymbol{\mathcal{M}(\theta_0,\theta_1,\theta_\infty)}$} \label{4.1}

 As before we assume that the parameters satisfy ``${\rm restr}$''.
 We recall that the modules in ${\bf S}(\theta_0,\theta_1,\theta_\infty)$ are assumed to be non-split. For special values of the $\{\theta_*\}$, the space $\mathcal{M}(\theta_0,\theta_1,\theta_\infty)$, constructed in Section~\ref{Section3}, contains split modules.
 {\it In such a case we delete the locus of the split modules and keep the notation $\mathcal{M}(\theta_0,\theta_1,\theta_\infty)$}.
 We compute these special values $\{\theta_*\}$ and the locus of the split modules.

 Let $M$ denote such a split module. Then $M=L\oplus L_2$, where $L$, $L_2$ are the only submodules of dimension 1.
 Both $L$ and $L_2$ correspond to saturated line bundles in $\mathcal{V}$ with degree $-1$.
 We may choose $e_2$ such that $L_2=\mathbb{C}(z)e_2$ and the corresponding line bundle is $O(-[\infty ])e_2\subset \mathcal{V}$.

 As above, the connection on $\mathcal{V}$ is represented by the differential operator
 \[
 D:=\frac{{\rm d}}{{\rm d}z}+\frac{1}{z(z-1)} \begin{pmatrix} a& 0\\ b& -a\end{pmatrix}
 \]
 with $a=a_0+a_1z+a_2z^2$, $b=b_1z+b_0\neq 0$ and as before
\[
a_0=\epsilon_0\frac{\theta_0}{2}, \qquad a_1=\epsilon_1\frac{\theta_1}{2}-\epsilon_0\frac{\theta_0}{2}-\epsilon_2\frac{t}{2},
\qquad a_2=\epsilon_2\frac{t}{2}.
\]

 Let $\mathcal{L}\subset \mathcal{V}$ denote the (saturated) line bundle corresponding to $L$. The line bundle
 $\mathcal{L}([\infty ])$, is saturated in $ \mathcal{V}([\infty])=O([\infty ])e_1+Oe_2$ and has degree zero. It is generated by an element $v=\alpha e_1+e_2$ with $\alpha =d_1z+d_0$. The assumption $D(v)$ is a multiple of $v$ leads to the equation
\[
\frac{{\rm d}}{{\rm d}z}(\alpha)+\frac{\alpha a}{z(z-1)}=\frac{\alpha^2b-\alpha a}{z(z-1)}.
\]
 This implies that $\alpha =d_1z$ or
 $\alpha =d_1(z-1)$ and, by multiplying $e_2$ by a constant, we may assume that $d_1=1$.
 We continue with the first case, which translates into the equations
 \[
 2a_2-b_1=0, \qquad 1+2a_1-b_0=0, \qquad -1+2a_0=0.
 \]
 The last equality implies $\epsilon_0\theta_0=1$ and $\epsilon_2\theta +\epsilon_1\theta_1=0$.
 Hence $\theta_0=\pm 1$ and $\theta=\pm \theta_1$. Moreover, $b_0=1+2a_1$ and $b_1=2a_2$.
 Thus this part of the split locus is isomorphic to $\operatorname{Spec}\big(\mathbb{C}\big[t,t^{-1}\big]\big)$.

 The same holds for the other part of the split locus corresponding to $\alpha =(z-1)$ and
 $\epsilon_1\theta_1=1$ and $\epsilon_2\theta -\epsilon_0\theta_0=0$. Hence $\theta_1=\pm 1$ and $\theta=\pm \theta_0$.

{\it A typical example} is $\theta_0=1$, $\epsilon_0=1$, ``generic'' $\theta_1$ and $\theta=\theta_1$.
The two data for submodules with degree $-1$ are $(\epsilon_0, \epsilon_1,\epsilon_2)=(1,1,-1)$ and
$(\epsilon_0, \epsilon_1,\epsilon_2)=(1,-1,1)$. The reducible locus in the original $\mathcal{M}(\theta_0,\theta_1,\theta_\infty)$ consists of two projective lines
$\mathbb{P}^1_{\mathbb{C}[t,t^{-1}]}\rightarrow \operatorname{Spec}\big(\mathbb{C}\big[t,t^{-1}\big]\big)$.
The intersection of these projective lines is one point, equal to the locus of the split modules. Removing the split locus results in two non-intersecting affine lines over the same ring.
Summarizing, the reducible locus of the new $\mathcal{M}(\theta_0,\theta_1,\theta_\infty)$ consists of two non-intersecting affine lines over~$\mathbb{C}\big[t,t^{-1}\big]$.

 We note that the special cases coincide with two of the three exceptional cases for the monodromy space in part~(v) on page~\pageref{remv}, where reducibility with both $\epsilon_2=1$ and $\epsilon_2=-1$ occur.

In accordance with what has been done for the monodromy space we introduce for all {\it three exceptional cases the moduli spaces of connections}
$\mathcal{M}(\theta_0,\theta_1,\theta_\infty, \epsilon_2\!=\!-1)$ and $\mathcal{M}(\theta_0,\theta_1,\theta_\infty, {\epsilon_2\!=\!1})$.
The first one is obtained by deleting the locus of the split objects and only admitting the reducible objects
with $\epsilon_2=-1$. The second space
is defined by replacing $\epsilon_2=-1$ by $\epsilon_2=1$.%
\label{Mrefined}

\begin{Remark} For the exceptional case $(\theta =0, -\epsilon_0 \theta_0+\epsilon_1\theta_1=-1 )$, the locus of reducible objects with
$\epsilon_2=-1$ coincides with the one for $\epsilon_2=1$. The same holds for the locus of the split objects. For these parameters, the space
 $\mathcal{M}(\theta_0,\theta_1,\theta_\infty, \epsilon_2=-1)$ is obtained by deleting only the locus of the split modules.
\end{Remark}

 \section{The parabolic structure}\label{Section5}

 {\it Parabolic structure for connections and another way to formulate Lemma~$\ref{2.1}$}.
 As before, we assume that the parameters satisfy ``${\rm restr}$''.
 Consider a regular singular differential module $N=(N,D)$ over $\mathbb{C}(\!(z)\!)$ with
 $D(f\cdot n)=z\frac{{\rm d}f}{{\rm d}z}\cdot n +f\cdot D(n)$, $\dim N=2$ and $\det N={\bf 1}$. An element
 $\eta \in \mathbb{C}$ will be called an {\it eigenvalue for $N$} if there is a $D$-invariant lattice $\Lambda$ for $N$
 and $\eta$ is an eigenvalue of $D$ on $\Lambda /z\Lambda$.

 Suppose that $\eta$ is an eigenvalue. Then $N$ is represented by a differential operator of the form
 $z\frac{{\rm d}}{{\rm d}z}+\big(\begin{smallmatrix} a& b\\ c& d\end{smallmatrix}\big)$ with $a,b,c,d\in
 \mathbb{C}[\![z]\!]$,
 $ \big(\begin{smallmatrix} a& b\\ c& d\end{smallmatrix}\big) \bmod (z)$ has eigenvalue $\eta$ and eigenvalue
 $-\eta +m$ for some integer~$m$. One easily concludes that the set of all eigenvalues is
 $(\eta +\mathbb{Z})\cup (-\eta +\mathbb{Z})$.

 A {\it parabolic structure} for $(N,\eta)$ is an {\it eigenline} $\mathbb{C}e$ with
 $e\in N,\ e\neq 0$ and $D(e)=\eta e$. If $(N,\eta)$ has a parabolic structure, then $\eta$ is an
 eigenvalue for $N$. Now there are the following possibilities.
\begin{enumerate}\itemsep=0pt
 \item[$({\rm i})$] $N$ has eigenvalue $\eta$ and $2\eta \not \in \mathbb{Z}$. Then
 $(\eta +\mathbb{Z})\cap (-\eta +\mathbb{Z})=\varnothing$ and $(N,\alpha)$ has a unique
 parabolic structure
 for every $\alpha \in (\eta +\mathbb{Z})\cup (-\eta +\mathbb{Z})$.
 Indeed, the condition $2\eta \not \in \mathbb{Z}$ implies that $N$ is a direct sum, i.e.,
 corresponds to an operator of the form
 $z\frac{{\rm d}}{{\rm d}z}+\big(\begin{smallmatrix} \eta& 0\\ 0& -\eta\end{smallmatrix}\big)$.

\item[$({\rm ii})$] $\eta =0$ is eigenvalue of $N$. Then one can represent $N$ by the operator
 $z\frac{{\rm d}}{{\rm d}z}+\big(\begin{smallmatrix} 0& a\\ 0& 0\end{smallmatrix}\big)$ with $a\in \{0,1\}$.
\begin{enumerate}\itemsep=0pt
\item[$({\rm a})$] If $a=0$, then $\ker D$ is a two-dimensional complex vector space and
 the parabolic structures for $\eta=0$ are the lines in $\ker D$. Thus the level structures are parametrized by
 $\mathbb{P}^1$ (over $\mathbb{C}$).

\item[$({\rm b})$] If $a=1$, then $\ker D$ has dimension 1 over $\mathbb{C}$. This is the only parabolic structure for
 $\eta=0$.
\end{enumerate}

\item[$({\rm iii})$] $\eta =\frac{m}{2}$ and $m\in \mathbb{Z}_{>0}$. It can be shown that $N$ is represented
 by an operator $z\frac{{\rm d}}{{\rm d}z}+\Big(\begin{smallmatrix} -\frac{m}{2}& az^m\\ 0& \frac{m}{2}\end{smallmatrix}\Big) $ with
 $a\in \{0,1\}$, say on a basis $e_1$, $e_2$ for $N$. For $(N,\frac{m}{2})$ there is a unique parabolic structure.
\begin{enumerate}\itemsep=0pt
\item[$({\rm a})$] If $a=0$, then the parabolic structures for $-\frac{m}{2}$ are the eigenlines
 $\mathbb{C}(c_1e_1+c_2z^{-m}e_2)$. The eigenlines are parametrized by $\mathbb{P}^1$.

\item[$({\rm b})$] If $a=1$, then the only parabolic structure for $-\frac{m}{2}$ is $\mathbb{C}e_1$.
\end{enumerate}
\end{enumerate}

 Recall that ${\bf S}(\theta_0,\theta_1,\theta_\infty)$ denotes the set of equivalence classes of differential modules
 with prescribed singularities and provided with a parabolic structure.

{\it The definition of $\mathcal{M}(\theta_0,\theta_1,\theta_\infty)$ for $\theta_0\in \mathbb{Z}$ and/or
 $\theta_1\in\mathbb{Z}$ with parabolic structure. {\em Example:} the case $\theta_0\in \mathbb{Z}$ and $\theta_1\not \in \mathbb{Z}$}.
 Consider the affine space $M(\theta_0,\theta_1,\theta_\infty)$ defined on
 page~\pageref{DefM} by variables $a_*$, $b_*$, $c_*$, $d_*$
 and many relations. Take the subspace of $\mathbb{P}^1\times M(\theta_0,\theta_1,\theta_\infty)$
 consisting of the pairs
 \[\left([y_1:y_2], \frac{{\rm d}}{{\rm d}z}+\frac{1}{z(z-1)}\cdot
 \left(\begin{matrix} a_0+a_1z+a_2z^2& c_0+c_1z+c_2z^2+c_3z^3 \\ b_0+b_1z& d_0+d_1z+d_2z^2 \end{matrix}\right) \right)
\]
with $\big(\begin{smallmatrix} -a_0 & -c_0\\ -b_0& -d_0\end{smallmatrix}\big) \big(\begin{smallmatrix} y_1\\ y_2\end{smallmatrix}\big) \equiv \frac{-|\theta_0|}{2}\big(\begin{smallmatrix} y_1\\ y_2\end{smallmatrix}\big) \bmod (z)$ and $b_0+b_1z\neq 0$.
 We observe that there is a unique eigenline for $-\frac{|\theta_0|}{2}$ which reduces modulo $(z)$ to the above
 element $[y_1:y_2]\in \mathbb{P}^1$. Therefore,
 the quotient of this space by the action of the automorphism group $G$ of the vector bundle $\mathcal{V}$ is the
 moduli space we are looking for.

 The above space can be cut into four open affine subsets given by inequalities
 ($y_1\neq 0$ or $y_2\neq 0$) and ($b_0\neq 0$ or $b_1\neq 0$). On each of these affine parts, taking the
 geometric quotient by~$G$ is obtained by normalizing ($y_1=1$ or $y_2=1$) and ($b_0=1$ or $b_1=1$).

 The above definition and construction holds in all cases with $\theta_0\in \mathbb{Z}$ and/or $\theta_1\in \mathbb{Z}$.
 One easily verifies that by gluing one obtains an irreducible smooth variety of dimension~3,
 denoted by $\mathcal{M}(\theta_0,\theta_1,\theta_\infty)$, such that the natural map
 $\mathcal{M}(\theta_0,\theta_1,\theta_\infty)(\mathbb{C})\rightarrow {\bf S}(\theta_0,\theta_1,\theta_\infty)$ is a bijection.

 For, say, $\theta_0\in 2\mathbb{Z}$, an eigenline for the connection at $z=0$ is also an eigenline for the local
 monodromy at $z=0$ for $s_1=2$. This holds for all cases with $\theta_0\in \mathbb{Z}$ and/or $\theta_1\in \mathbb{Z}$. Thus we have, for all $\theta_0$, $\theta_1$, $\theta_\infty$, an analytic morphism ${\rm RH}\colon \mathcal{M}(\theta_0,\theta_1,\theta_\infty)\rightarrow \mathcal{R}^{\rm geom}(\theta_0,\theta_1,\theta_\infty)$. Using Proposition~\ref{2.2}, one easily verifies that the proof of Theorem~\ref{2.4} extends to a proof of
 Theorem~\ref{2.5}.

 We note that for exceptional parameters, the two families of spaces
 $\mathcal{M}(\theta_0,\theta_1,\theta_\infty,\epsilon_2=\pm 1)$ are introduced and studied in Section~\ref{4.1}.

 \begin{Theorem}\label{2.5} Assume that the parameters satisfy ``${\rm restr}$'', i.e.,
 $-1/2-\epsilon_2\frac{\theta}{2}+\epsilon_1\frac{\theta_1}{2}-\epsilon_0\frac{\theta_0}{2}\in \mathbb{Z}_{<-1}$ is not possible
 for $(\epsilon_0,\epsilon_1,\epsilon_2)\in \{\pm 1\}^3$.
 \begin{enumerate}\itemsep=0pt
 \item[$(1)$] If the parameters are not exceptional, then the extended Riemann--Hilbert map
\[{\rm RH}^{+}\colon\ \mathcal{M}(\theta_0,\theta_1,\theta_\infty)\rightarrow
\mathcal{R}^{\rm geom}(\theta_0,\theta_1,\theta_\infty)\times \mathbb{C}\]
is an analytic isomorphism.
\item[$(2)$] If the parameters are exceptional, then the above holds after adding the data $\epsilon_2=-1$ or $\epsilon_2=1$
to both spaces.
\end{enumerate}
\end{Theorem}

 \section{The Lax pair and the fifth Painlev\'e equation}\label{Section6}

The explicit connection $\frac{{\rm d}}{{\rm d}z}+A$ with fixed $\theta_0$, $\theta_1$, $\theta$ given on page~\pageref{Dmatrix}
(below Proposition~\ref{2.3}), can be completed
 to a Lax pair by using $\frac{{\rm d}}{{\rm d}t}+B$, where the entries of $B$ are polynomials in~$z$ of degrees $\leq 2$
over a field of meromorphic functions in $t$. The relation $\big[\frac{{\rm d}}{{\rm d}z}+A,\frac{{\rm d}}{{\rm d}t}+B\big]=0$ determines $A$ and $B$
completely. The zero $q=-b_0$ of the entry $A[2,1]$ (2nd row, 1st column) satisfies the fifth Painlev\'e equation corresponding to
the parameters $\theta_0$, $\theta_1$, $\theta_\infty +1$ (the $+1$ comes from the shift we made in the local data at
$z=\infty$). This equation reads
\begin{small}
\[q'' = \frac{1}{2}\left(\frac{1}{q}+\frac{1}{q - 1}\right)(q')^2-\frac{q'}{t}+\frac{q(q - 1)}{t}\theta-\frac{q - 1}{2qt^2}\theta_0^2
+\frac{q}{2t^2(q - 1)}\theta_1^2 -\frac{(2q - 1)(q - 1)q}{2}.\]
\end{small}
Substituting $q=\frac{y}{y-1}$ produces for $y$ the standard Painlev\'e equation
\[ {\rm P}_5\left(\alpha=\frac{\theta_1^2}{2}, \beta =-\frac{\theta_0^2}{2}, \gamma =-\theta, \delta=-\frac{1}{2}\right) \quad\mbox{and}\quad
\theta =1+\theta_\infty . \]

The formula here seems to correct the one for ${\rm P}_5$ in \cite[equation~(4.3)]{vdP-Sa}.

 \section{Connections for the Lax pair of Noumi--Yamada}\label{Section7}

For certain tuples $( \theta_0,\theta_1,\theta_\infty)$ and $( \theta^*_0,\theta^*_1,\theta^*_\infty)$,
there are birational maps $\mathcal{M}(\theta_0,\theta_1,\theta_\infty)\;\cdots\longrightarrow \mathcal{M}(\theta^*_0,\theta^*_1,\theta^*_\infty)$. These birational maps are the so-called B\"acklund transformations.
Their action on solutions can be made explicit. The B\"acklund transformations form the
extended affine Weyl group of $A_3^{(1)}$, see \cite{O2,N-Y}.

 From Theorem~\ref{2.5}, one can deduce that for instance shifts of $\theta_0$, $\theta_1$, $\theta_\infty$ over integers are B\"acklund transformations. However in our construction of $\mathcal{M}(\theta_0,\theta_1,\theta_\infty)$
 {\it some B\"acklund transformations are not easily visible}. Noumi and Yamada \cite{N-Y} (and others) studied a different Lax pair for ${\rm P}_5$ from which one can read off all B\"acklund transformations. Here we introduce a~natural moduli space of connections which leads to this Lax pair.

\subsection[Moduli spaces for q\_0=z\^{}\{1/2\}+tz\^{}\{1/4\} with t in C\^{}*]{Moduli spaces for $\boldsymbol{q_0=z^{1/2}+tz^{1/4}}$ with $\boldsymbol{t\in \mathbb{C}^*}$}
Consider the set $\mathbf{S_4}$ of isomorphism classes of differential modules $M$ over $\mathbb{C}(z)$ given by
\begin{enumerate}\itemsep=0pt
\item[$({\rm i})$] $\dim M=4$ and $\Lambda ^4M$ is the trivial differential module and the only singularities are $0$, $\infty$,
\item[$({\rm ii})$] $z=0$ is a regular singularity,
\item[$({\rm iii})$] $z=\infty$ is irregular and has eigenvalues $q_0$ and its conjugates $q_1$, $q_2$, $q_3$.
\end{enumerate}
 We want to convert the set $\mathbf{S_4}$ into a family of connections $(\nabla ,\mathcal{V})$ on a vector bundle $\mathcal{V}$ on $\mathbb{P}^1$ of rank 4 with singularities at $z=0$ and $z=\infty$, such that
$\nabla\colon \mathcal{V}\rightarrow \Omega _{\mathbb{P}^1}([0]+2[\infty ])\otimes \mathcal{V}$, or equivalently
$\delta=\nabla_{z\frac{{\rm d}}{{\rm d}z}}\colon \mathcal{V}\rightarrow \mathcal{V}([\infty ])$, has a
(chosen) prescribed
singularity at $z=\infty$.

A possible {\it prescribed singularity at $z=\infty$} is obtained as follows.
 For the fields $\mathbb{C}\big(z^{1/4}\big)$ and $\mathbb{C}\big(\!\big(z^{-1/4}\big)\!\big)$, we write $\sigma$ for
the automorphism given by $\sigma \big(z^{1/4}\big)={\rm i}z^{1/4}$. Let $N$ denote the free
$\mathbb{C}[\![z^{-1/4}]\!]$ module with basis $r_0$, $r_1$, $r_2$, $r_3$ provided with a semi-linear action of $\sigma$
such that $\sigma$ satisfies $r_0\mapsto r_1\mapsto r_2\mapsto r_3\mapsto r_0$.
 Let $N_0\subset N^{\langle \sigma\rangle}$ be the $\mathbb{C}[\![z^{-1}]\!]$-module with free basis
$s_0=r_0+r_1+r_2+r_3$, $s_1=z^{1/4}(r_0+{\rm i}r_1-r_2-{\rm i}r_3)$, $s_2=z^{1/2}(r_0-r_1+r_2-r_3)$, $s_3=z^{3/4}(r_0-{\rm i}r_1-r_2+{\rm i}r_3)$. We note that $\mathbb{C}\big(\!\big(z^{-1}\big)\!\big)\otimes N_0=\mathbb{C}\big(\!\big(z^{-1}\big)\!\big)\otimes N^{\langle \sigma\rangle}$.

One defines $\delta_0\colon N\rightarrow N\otimes \mathbb{C}\big(\!\big(z^{-1/4}\big)\!\big)$
by $\delta_0(r_j)=\big(q_j-\frac{3}{8}\big)r_j$ for $j=0,1,2,3$
and $\delta_0(fn)=z\frac{{\rm d}f}{{\rm d}z}n+f\delta_0(n)$ for $f\in \mathbb{C}[\![z^{-1/4}]\!]$, $n\in N$.
The induced operator $N_0\rightarrow N_0\otimes \mathbb{C}\big(\!\big(z^{-1}\big)\!\big)$ is also denoted by $\delta_0$. The term $-\frac{3}{8}$ is introduced in order that the operator on $N_0$
has trace zero. We note that $(N_0,\delta_0)$ is irreducible and therefore its group of automorphisms is
$\mathbb{C}^*$. The chosen prescribed singularity at $z=\infty$ is the formal local connection $(N_0,\delta_0)$.
 On the basis $\{s_0,s_1,s_2,s_3\}$ of $N_0$, this is represented by
\[
L_\infty:=z\frac{{\rm d}}{{\rm d}z}+{ \left(\begin{matrix} -\frac{3}{8}& 0& z& \frac{t}{4}z\vspace{1mm}\\ \frac{t}{4}& -\frac{1}{8} & 0 &z \vspace{1mm}\\ 1&\frac{t}{4} &\frac{1}{8} & 0\vspace{1mm}\\ 0& 1&\frac{t}{4} &\frac{3}{8} \end{matrix}\right)}.
\]

\begin{Remark}\label{Rem7.1}
The choice $(N_0,\delta_0)$ and corresponding $L_\infty$ used here is not unique. A canonical choice will
be described in Remark~\ref{canonicalchoice}.
\end{Remark}

We consider {\it the moduli functor} for the data:
``The isomorphy classes of the connections $\delta \colon \mathcal{V}\rightarrow \mathcal{V}([\infty ])$ on the {\it free bundle} $\mathcal{V}$ such that the induced $\delta\colon \widehat{\mathcal{V}}_\infty \rightarrow \widehat{\mathcal{V}}_\infty\otimes \mathbb{C}\big(\!\big(z^{-1}\big)\!\big)$ is isomorphic to $(N_0,\delta_0)$''.

This isomorphism is unique up to $\mathbb{C}^*$. The connections $\delta$ (modulo isomorphisms) can be represented by the operators $L:=z\frac{{\rm d}}{{\rm d}z}+A_0+A_1z$ (modulo equivalence) for which there exists a
$U\in {\rm GL}_4\big(\mathbb{C}[\![z^{-1}]\!]\big)$ with $ULU^{-1}=L_\infty$.

The above is an example for the theory developed in~\cite[Section~12 (especially Section~12.5)]{vdP-S}. It is shown there that a fine moduli space exists by giving an explicit construction for the general ramified case. {\it For the present case, this moduli space is denoted by $\mathcal{M}_4$} and we will produce explicit formulas.

Consider the map $\eta\colon \mathcal{M}_4(\mathbb{C})\rightarrow {\mathbf{S_4}}$ associating to a~connection
its generic fibre (which is a~differential module over $\mathbb{C}(z)$ belonging to $\mathbf{S_4}$). From the observation
that $(N_0,\delta_0)$ has automorphism group $\mathbb{C}^*$, it follows that $\eta$ is injective.

For $M\in {\mathbf{S_4}}$, an invariant lattice $\Lambda_\infty$ at $z=\infty$ is defined as follows:
$\Lambda_\infty$ is a free $\mathbb{C}[\![z^{-1}]\!]$-submodule of $\mathbb{C}\big(\!\big(z^{-1}\big)\!\big)\otimes M$ of rank 4, such that $\delta (\Lambda_\infty)\subset z\Lambda_\infty $.
Similarly, an invariant lattice~$\Lambda_0$ for~$M$ at $z=0$ is by definition
a free $\mathbb{C}[\![z]\!]$-submodule of $\mathbb{C}(\!(z)\!)\otimes M$ of rank 4 such that
$\delta(\Lambda_0)\subset \Lambda_0$.

 From the definition of $\mathbf{S_4}$, it follows
 that any $M\in \mathbf{S_4}$ has a {\it unique} invariant lattice $\Lambda_\infty\subset \mathbb{C}\big(\!\big(z^{-1}\big)\!\big)\otimes M$ at $z=\infty$ which is isomorphic to $(N_0,\delta_0)$. Further, one chooses an invariant lattice
$\Lambda_0\subset \mathbb{C}(\!(z)\!)\otimes M$. These lattices determine a connection $(\nabla, \mathcal{V})$
with generic fibre $M$. If $\mathcal{V}$ happens to be the free vector bundle, then $M$ lies in the image of $\eta$.

The regular singular differential module $\mathbb{C}(\!(z)\!)\otimes M$ can be represented by a matrix differential operator $z\frac{{\rm d}}{{\rm d}z}+A$, such that the eigenvalues of the constant matrix $A$ are distinct modulo integers.
If $A$ is diagonalizable, then $\mathbb{C}(\!(z)\!)\otimes M$ contains many invariant lattices and, for a suitable choice of the lattice, the corresponding vector bundle $\mathcal{V}$ is free (compare the method of \cite[Section~6.5]{vdP-S}).
This shows that the image of~$\eta$ contains what could be called a ``Zariski open, dense subset of~$\mathbf{S_4}$''. However, the map~$\eta$ is probably not surjective.

The parameter space $\mathcal{P}^+_4\cong \mathbb{C}^3$ is the space of the polynomials of the form $T^4+c_2T^2+c_1T+c_0$. The morphism $\mathcal{M}_4\rightarrow \mathcal{P}_4^+$ maps a connection
$z\frac{{\rm d}}{{\rm d}z}+A_0+zA_1\in \mathcal{M}_4(\mathbb{C})$ to the characteristic polynomial of $A_0$.

Let $\mathcal{R}_4$ denote the moduli space for the analytic data (i.e., the topological monodromy and the Stokes matrices) of the set ${\mathbf{S_4}}$. The Riemann--Hilbert morphism is the obvious map ${{\rm RH}\colon \mathcal{M}_4\rightarrow \mathcal{R}_4}$. The map ${\rm RH}$ forgets $t$ and the choice of a logarithm of the topological monodromy at $z=0$. Thus
the fibres of ${\rm RH}$ have dimension 1 and so $\dim \mathcal{M}_4=1+\dim \mathcal{R}_4$.
The connected components of the fibres are parametrized by $t$.

 We note that for every integer $n\geq 3$, one can define, in a similar way, a family of differential modules
$\bf S_n$, given as follows: $\dim M=n$; $\Lambda^nM$ is the trivial differential module; $z=0$ is regular singular and the (generalized) eigenvalues at $z=\infty$ are $z^{2/n}+tz^{1/n}$ and its conjugates. The corresponding moduli spaces $\mathcal{M}_n$, $\mathcal{R}_n$ produce the Lax pairs introduced by Noumi, Yamada et al. For $n=3$ this Lax pair yields~${\rm P}_4$. This is exploited in \cite{vdP-T}. For $n=4$, the Lax pair produces~${\rm P}_5$ and its B\"acklund transformations. We will give the details for $n=4$ and refer to a~future paper for the general case.

\subsubsection[R\_4 and the fibres of R\_4 rightarrow P\_4]{$\boldsymbol{\mathcal{R}_4}$ and the fibres of $\boldsymbol{\mathcal{R}_4\rightarrow \mathcal{P}_4}$}
 The monodromy of a differential module $M\in {\mathbf{S_4}}$ consists of the topological monodromies
 ${\rm mon}_0$ and ${\rm mon}_\infty$ at the points $0$ and $\infty$. Since there are no other singularities
 one has ${\rm mon}_0\cdot {\rm mon}_\infty =1$. The monodromy at $\infty$ is given by the Stokes matrices and
 the formal monodromy $\gamma$. The monodromy identity \cite{vdP-S} states that ${\rm mon}_\infty$
 is equal to their product.

Now we compute the Stokes data at $z=\infty$ and refer to \cite{vdP-S} for more details.
 The (generalized) eigenvalues at $z=\infty$ are:
\begin{gather*}
q_0 = z^{1/2}+tz^{1/4},\qquad q_1 = -z^{1/2}+{\rm i}tz^{1/4},\qquad q_2=z^{1/2}-tz^{1/4},\qquad q_3=-z^{1/2}-{\rm i}tz^{1/4}.
\end{gather*}
 The differences $q_0-q_1$, $q_0-q_3$, $q_2-q_1$, $q_2-q_3$ have the form $2z^{1/2}+\cdots$, and
 $q_0-q_2=2tz^{1/4}$ and $q_1-q_3=2{\rm i}tz^{1/4}$. There is one singular direction in $[0,1)$ for the terms
 $\pm 2z^{1/2}$ as well as for each of $\pm 2tz^{1/4}$, $\pm 2{\rm i}z^{1/4}$.
 This leads to the monodromy identity
 \[{\rm mon}_\infty= \left(\begin{matrix} & & &-1 \\ 1& & & \\ & 1& & \\ & &1 & \end{matrix}\right) \left(\begin{matrix} 1 & & & \\ & 1& &y \\ & &1 & \\ & & & 1\end{matrix}\right)
 \left(\begin{matrix} 1 & & & \\ x_1&1 &x_2 & \\ & & 1& \\ x_3& &x_4 & 1\end{matrix}\right). \]
 The basis for the above matrices is unique up to multiplying every basis vector by the same scalar. Therefore, $\mathcal{R}_4\cong \mathbb{C}^5$.

 The parameter space $\mathcal{P}_4$ for $\mathcal{R}_4$ is given
 by the characteristic polynomial of ${\rm mon}_\infty$, which is written as
 $T^4+p_3T^3+p_2T^2+p_1T+1$. Thus $\mathcal{P}_4\cong \mathbb{C}^3$.
 One considers the commutative diagram
 \[\begin{array}{@{}ccc} \mathcal{M}_4& \stackrel{{\rm RH}}{\rightarrow} &\mathcal{R}_4\\
 \downarrow & & \downarrow \\
 \mathcal{P}_4^+ & \stackrel{{\rm rh}}{\rightarrow}& \mathcal{P}_4, \end{array}
 \]
 where ${\rm rh}$ maps the characteristic polynomial of a local differential operator
 $z\frac{{\rm d}}{{\rm d}z}+A$ at $z=0$ to the characteristic polynomial of ${\rm mon}_\infty ={\rm mon}_0^{-1}$;
 explicitly, $\prod_{j=1}^4(T-\lambda _j) \mapsto \prod_{j=1}^4\big(T-{\rm e}^{2\pi {\rm i} \lambda _j}\big)$.

 The fibre of $\mathcal{R}_4\rightarrow \mathcal{P}_4$ above $(p_1,p_2,p_3)$ is defined by an obvious set
 of equations in the variables $x_1$, $x_2$, $x_3$, $x_4$, $y$. After elimination of the two variables $x_1$, $x_2$, there remains
 one cubic equation in the variables $y$, $x_3$, $x_4$ and parameters $p_1$, $p_2$, $p_3$. After a linear change of $y$, $x_3$, $x_4$
the equation has the form
 \[v_1v_2v_3+*v_1^2+*v_2^2+*v_1+*v_2+*v_3+*=0\] for suitable affine expressions $*$'s in the parameters
 $p_1$, $p_2$, $p_3$. The cubic equation for the monodromy of ${\rm P}_5$ has the same features. Comparing with the table of cubic surfaces associated to Painlev\'e equations in \cite[pp.~2635--2636]{vdP-Sa} and \cite{C-M-R,PR23}, this
is an indication that the (yet to be determined) Lax pair produces ${\rm P}_5$.

\subsubsection[A matrix differential operator for M\_4, Lax pairs and P\_5]{A matrix differential operator for $\boldsymbol{\mathcal{M}_4}$, Lax pairs and $\boldsymbol{{\rm P}_5}$}
We want to represent a ``general'' element $M\in \mathcal{M}_4$ by a matrix differential operator
$z\frac{{\rm d}}{{\rm d}z}+A_0+zA_1$. It turns out that a direct computation is rather difficult. Instead,
we adapt the computation of the formal local connection $(N_0,\delta_0)$ and $L_\infty$, to
the global case. Then $M$ is replaced by $N:=\mathbb{C}\big(z^{1/4}\big)\otimes M$ and we
use ``equivariant formules'' on $N$.
Let $\gamma$ denote the automorphism of $\mathbb{C}\big(z^{1/4}\big)$ over $\mathbb{C}(z)$, defined by
$\gamma \big(z^{1/4}\big)={\rm i}z^{1/4}$. Let $\sigma\colon N\rightarrow N$ be given by
$\sigma (f\otimes m)=\gamma(f)\otimes m$ for $f\in \mathbb{C}\big(z^{1/4}\big)$ and $m\in M$.
This new differential module $N$ has no ramification
at $z^{1/4}=\infty$ and has Katz invariant 2. Let $D$ denote the differential operator on~$N$.
One considers a basis $e_0$, $e_1$, $e_2$, $e_3$ of $N$ such that $\sigma$ acts by
$e_0\mapsto e_1\mapsto e_2\mapsto e_3\mapsto e_0$. We try to construct $D$ by a formula
 \[D(e_0)=\bigg(z^{1/2}+\frac{t}{4} z^{1/4}-3/8\bigg) e_0 +\big(a_1+b_1z^{1/4}\big) e_1+a_2 e_2+\big(a_3+b_3z^{1/4}\big) e_3. \]
 This determines $D$ because $\sigma\circ D=D\circ \sigma$. The term $-3/8$ is introduced
 (as before) in order to obtain a matrix for $D$ with trace zero.

 The formula for $D$ follows from the construction of the moduli space for the
 ramified case \cite[Section~12.5]{vdP-S}. The matrix of $D$ on the basis of invariants
\begin{gather*}
 B_0=e_0+e_1+e_2+e_3,\\
 B_1=z^{1/4}\big(e_0+{\rm i}e_1+{\rm i}^2e_2+{\rm i}^3e_3\big),\\
 B_2=z^{1/2}(e_0-e_1+e_2-e_3),\\
 B_3=z^{3/4}(e_0-{\rm i}e_1-e_2+{\rm i}e_3)
 \end{gather*}
 is
\[\left(\begin{matrix} -\frac{3}{8}+a_1+a_2+a_3 & 0 & z & z\big(\frac{t}{4}+b_1+b_3\big) \vspace{1mm}\\
 \frac{t}{4}-{\rm i}b_1+{\rm i}b_3 & -\frac{1}{8}-{\rm i}a_1-a_2+{\rm i}a_3 & 0 &z \vspace{1mm}\\
 1 &\frac{t}{4}-b_1-b_3 &\frac{1}{8}-a_1+a_2-a_3 & 0 \vspace{1mm}\\
 0 & 1 & \frac{t}{4}+{\rm i}b_1-{\rm i}b_3 &\frac{3}{8}+{\rm i}a_1-a_2-{\rm i}a_3
 \end{matrix}\right) \]
 and $D$ is equal to the differential operator
\[
z\frac{{\rm d}}{{\rm d}z}+\left(\begin{matrix} \epsilon_0&0&z&zf_0\\
 f_1&\epsilon_1&0&z\\ 1&f_2&\epsilon_2& 0\\ 0&1&f_3&\epsilon_3 \end{matrix}\right)
 \]
 with $\sum \epsilon _j=0$, $f_0+f_2=f_1+f_3=\frac{t}{2}$. The $\epsilon_0,\dots ,\epsilon_3$ are parameters.

 The operator $D$ is completed to a Lax pair by the differential operator $E$ with respect to $\frac{{\rm d}}{{\rm d}t}$. This
 operator, written on the basis $e_0$, $e_1$, $e_2$, $e_3$ is $\sigma$-invariant and has the form
 $E(e_0)=z^{1/4}e_0+\sum _{j=1}^3h_je_j$ for suitable functions $h_1$, $h_2$, $h_3$ of $t$. On the basis
 $B_0$, $B_1$, $B_2$, $B_3$ one obtains
\[
E:=\frac{{\rm d}}{{\rm d}t}+ \left(\begin{matrix} g_0& 0& 0&z\\ 1&g_1&0&0\\0&1&g_2&0\\0&0&1&g_3 \end{matrix}\right)
\]
 with $\sum g_j=0$. The assumption that $E$ commutes with $D$ produces a~nonlinear system of first order differential equations for
 $f_0$, $f_1$, $f_2$, $f_3$ as functions of $t$.

 These formulas are almost identical to those derived by Noumi--Yamada. After eliminating~$f_2$,~$f_3$
 a combination of the resulting differential equations for $f_0$ and $f_1$ leads to the standard~${\rm P}_5$ equation, see \cite{N-Y,S-H-C} for details. The differential equations for $f_0$, $f_1$ are
 \begin{gather*}
 2tf'_0 = 8f_0^2f_1 - 2f_0^2t - 4f_0f_1t + f_0t^2 - 4\epsilon_0f_0 + 2\epsilon_0t + 4\epsilon_1f_0 -
 4\epsilon_2f_0 \\
 \hphantom{2tf'_0 =}{}
 + 4\epsilon_3f_0 - 2\epsilon_3t - 2f_0 +2t,\\
 2tf'_1 = -8f_0f_1^2 + 4f_0f_1t + 2f_1^2t - f_1t^2 + 4\epsilon_0f_1 - 2\epsilon_0t - 4\epsilon_1f_1
 + 2\epsilon_1t \\
 \hphantom{2tf'_1 =}{}
 + 4\epsilon_2f_1 - 4\epsilon_3f_1 + 2f_1.
 \end{gather*}
 In fact, the equations coincide with the ones in \cite{S-H-C} for $f_0$, $f_1$, $f_2$, $f_3$. The symmetry of these equations with respect to the cyclic permutation $f_0\mapsto f_1\mapsto f_2\mapsto f_3\mapsto f_0$ produces an element~$\pi$ of order 4 in the group of B\"acklund transformations. This is the ``missing element'' needed to obtain the extended affine Weyl group of $A_3^{(1)}$.

An {\it alternative computation for the Lax pair} is based on observing that isomonodromy is given by
 $DE(e_0)=ED(e_0)$ and that $a_1$, $a_2$, $a_3$ are independent of $t$. A straightforward computation then yields
 \begin{gather*}
 4t\cdot \frac{{\rm d}b_1}{{\rm d}t} = -16{\rm i}b_1^2b_3 + t^2b_1{\rm i} + 16{\rm i}b_3^3 - 4{\rm i}a_1t + 4a_1t - 32a_2b_3,\\
 4t\cdot \frac{{\rm d}b_3}{{\rm d}t} = -16{\rm i}b_1^3 + 16{\rm i}b_1b_3^2 - t^2b_3{\rm i} + 4{\rm i}a_3t - 32a_2b_1 + 4a_3t,\\
 h_1 = -b_1{\rm i}/2 + b_1/2, \qquad t\cdot h_2 = 2b_1^2{\rm i} - 2b_3^2{\rm i} + 4a_2,\qquad h_3 = b_3{\rm i}/2 + b_3/2.
 \end{gather*}
These equations for $b_1$, $b_3$ form a Hamiltonian system
 \begin{gather*}
 t\frac{{\rm d}b_1}{{\rm d}t}=\frac{\partial H}{\partial b_3},\qquad t\frac{{\rm d}b_3}{{\rm d}t}=-\frac{\partial H }{\partial b_1},\\
H=\frac{{\rm i}\big(b_1^2-b_3^2\big)^2}{t}+\frac{{\rm i}tb_1b_3}{4}+\frac{4a_2\big(b_1^2-b_3^2\big)}{t}-(1+{\rm i})a_3b_1+(1-{\rm i})a_1b_3.
\end{gather*}

The above Hamiltonian $H$ is, up to a change of variables, the same as Okamoto's polynomial Hamiltonian for ${\rm P}_5$, see \cite[p.~265]{O1}.

\begin{Remark}\label{canonicalchoice}{\rm
Here we describe a more canonical choice for the invariant lattice (see page~\pageref{Rem7.1}).
Let $N$ be the free $\mathbb{C}[\![z^{-1/4}]\!]$-module with basis $e_0$, $e_1$, $e_2$, $e_3$.
The action of $\sigma$ is given by
$e_0\mapsto e_1\mapsto e_2\mapsto e_3 \mapsto e_0$ and $\sigma \big(z^{1/4}\big)={\rm i}z^{1/4}$. The operator $\delta_0$
on $N$ commutes with $\sigma$ and is determined by $\delta_0(e_0)=\big(z^{1/2}+\frac{t}{4}z^{1/4}+\frac{3}{8}\big)e_0$.
Choose for the lattice $\tilde{N}_0$ at $z=\infty$ the module~$N^\sigma$ of the $\sigma$-invariant elements of $N$. It has a basis
\begin{gather*}
B_0=e_0+e_1+e_2+e_3,\\
\tilde{B}_1=z^{-3/4}(e_0+{\rm i}e_1-e_2-{\rm i}e_3),\\
\tilde{B}_2=z^{-1/2}(e_0-e_1+e_2-e_3),\\
\tilde{B}_3=z^{-1/4}(e_0-{\rm i}e_1-e_2+{\rm i}e_3).
\end{gather*}
 The global operator $D$, commuting with $\sigma$, is given by
 \[D(e_0)=\bigg(z^{1/2}+\frac{t}{4}z^{1/4}+\frac{3}{8}\bigg)e_0+\big(a_1+b_1z^{1/4}\big)e_1+(a_2)e_2+\big(a_3+b_3z^{1/4}\big)e_3.\]
 The matrix of $D$ with respect to the basis $B_0$, $\tilde{B}_1$, $\tilde{B}_2$, $\tilde{B}_3$ is
 \[ \left(\begin{matrix} \frac{3}{8}+a_1+a_2+a_3& 0&1 &b_1+b_3+\frac{t}{4} \vspace{1mm}\\
 z\big(\frac{t}{4}-{\rm i}b_1+{\rm i}b_3\big) &-\frac{3}{8}-{\rm i}a_1-a_2+{\rm i}a_3 & 0&z \vspace{1mm}\\
 z&\frac{t}{4}-b_1-b_3 &-\frac{1}{8}-a_1+a_2-a_3 & 0\vspace{1mm}\\
 0& 1&\frac{t}{4}+{\rm i}b_1-{\rm i}b_3 &\frac{1}{8}+{\rm i}a_1-a_2-{\rm i}a_3 \end{matrix}\right).\]
 This can be rewritten as
 \[
 \left(\begin{matrix}\epsilon_0 &0 &1 &g_0 \\ zf_0& \epsilon_1& 0& z\\ z&g_1 &\epsilon_2 &0
 \\ 0&1 &f_1 &\epsilon_3 \end{matrix}\right)
 \] with entries satisfying $f_0+f_1=g_0+g_1=\frac{t}{2}$ and $\epsilon_0+\cdots +\epsilon_3=0$.

Using the Lax pair conditions, elimination of $f_1$, $g_1$, $\epsilon_3$ and writing $f=f_0$, $g=g_0$ a system of differential equations essentially the same as the one in~\cite{S-H-C} is obtained
 \begin{gather*}
 2tf'= -8f^2g+2f^2t+ 4fgt-ft^2+8\epsilon_0f-2\epsilon_0t+2\epsilon_1t+8\epsilon_2f-2f+2t,\\
 2tg'= 8fg^2-4fgt-2g^2t +gt^2- 8\epsilon_0g+4\epsilon_0t+2\epsilon_1t- 8\epsilon_2g+2\epsilon_2t+2g.
 \end{gather*}

}\end{Remark}

\subsection*{Acknowledgements}

We are indebted to the referees of earlier versions of this manuscript for their careful reading and for the useful
suggestions and questions we obtained from them. One referee suggested that the construction of our moduli spaces could lead to
nonseparated spaces. This problem has finally, after additional comments from the same referee, been solved by a small restriction on the parameters, a ``natural'' restriction on the set of reducible differential modules in the sets ${\bf S}(\theta_0,\theta_1,\theta_\infty)$ and corresponding restrictions for the spaces $\mathcal{R}^{\rm geom}(\theta_0,\theta_1,\theta_\infty)$ and $\mathcal{M}(\theta_0,\theta_1,\theta_\infty)$.

\pdfbookmark[1]{References}{ref}
\LastPageEnding


\begin{thebibliography}{999}
\footnotesize\itemsep=0pt

\bibitem{A-P-T}
Acosta-Hum\'anez P.B., van~der Put M., Top J., Isomonodromy for the degenerate
 fifth {P}ainlev\'e equation, \href{https://doi.org/10.3842/SIGMA.2017.029}{\textit{SIGMA}} \textbf{13} (2017), 029,
 14~pages, \href{https://arxiv.org/abs/1612.03674}{arXiv:1612.03674}.

\bibitem{Boalch}
Boalch P., Symplectic manifolds and isomonodromic deformations, \href{https://doi.org/10.1006/aima.2001.1998}{\textit{Adv.
 Math.}} \textbf{163} (2001), 137--205, \href{https://arxiv.org/abs/2002.00052}{arXiv:2002.00052}.

\bibitem{Boalch2}
Boalch P., Quasi-{H}amiltonian geometry of meromorphic connections,
 \href{https://doi.org/10.1215/S0012-7094-07-13924-3}{\textit{Duke Math.~J.}} \textbf{139} (2007), 369--405,
 \href{https://arxiv.org/abs/math.DG/0203161}{arXiv:math.DG/0203161}.

\bibitem{Boalch3}
Boalch P., Wild character varieties, meromorphic {H}itchin systems and {D}ynkin
 diagrams, in Geometry and {P}hysics. {V}ol.~{II}, Oxford University Press,
 Oxford, 2018, 433--454, \href{https://arxiv.org/abs/1703.10376}{arXiv:1703.10376}.

\bibitem{C-M-R}
Chekhov L.O., Mazzocco M., Rubtsov V.N., Painlev\'e monodromy manifolds,
 decorated character varieties, and cluster algebras, \href{https://doi.org/10.1093/imrn/rnw219}{\textit{Int. Math. Res.
 Not.}} \textbf{2017} (2017), 7639--7691, \href{https://arxiv.org/abs/1511.03851}{arXiv:1511.03851}.

\bibitem{C2}
Clarkson P.A., Painlev\'e equations~-- nonlinear special functions,
 \href{https://doi.org/10.1016/S0377-0427(02)00589-7}{\textit{J.~Comput. Appl. Math.}} \textbf{153} (2003), 127--140.

\bibitem{C}
Clarkson P.A., Special polynomials associated with rational solutions of the
 fifth {P}ainlev\'e equation, \href{https://doi.org/10.1016/j.cam.2004.04.015}{\textit{J.~Comput. Appl. Math.}} \textbf{178}
 (2005), 111--129.

\bibitem{C3}
Clarkson P.A., Painlev\'e equations~-- nonlinear special functions, in
 Orthogonal {P}olynomials and {S}pecial {F}unctions, \textit{Lecture Notes in
 Math.}, Vol. 1883, \href{https://doi.org/10.1007/978-3-540-36716-1_7}{Springer}, Berlin, 2006, 331--411.

\bibitem{G}
Gromak V.I., On the transcendence of the {P}ainlev\'e equations,
 \textit{Differ. Equ.} \textbf{32} (1996), 156--162.

\bibitem{Hart}
Hartshorne R., Algebraic geometry, \textit{Grad. Texts Math.}, Vol.~52,
 \href{https://doi.org/10.1007/978-1-4757-3849-0}{Springer}, New York, 1977.

\bibitem{JM}
Jimbo M., Miwa T., Monodromy preserving deformation of linear ordinary
 differential equations with rational coefficients.~{II}, \href{https://doi.org/10.1016/0167-2789(81)90021-X}{\textit{Phys.~D}}
 \textbf{2} (1981), 407--448.

\bibitem{J-M-U}
Jimbo M., Miwa T., Ueno K., Monodromy preserving deformation of linear ordinary
 differential equations with rational coefficients.~{I}. {G}eneral theory and
 {$\tau$}-function, \href{https://doi.org/10.1016/0167-2789(81)90013-0}{\textit{Phys.~D}} \textbf{2} (1981), 306--352.

\bibitem{Mi}
Miwa T., Painlev\'e property of monodromy preserving deformation equations and
 the analyticity of {$\tau$} functions, \href{https://doi.org/10.2977/prims/1195185270}{\textit{Publ. Res. Inst. Math. Sci.}}
 \textbf{17} (1981), 703--721.

\bibitem{N-Y}
Noumi M., Yamada Y., Higher order {P}ainlev\'e equations of type {$A^{(1)}_l$},
 \textit{Funkcial. Ekvac.} \textbf{41} (1998), 483--503,
 \href{https://arxiv.org/abs/math.QA/9808003}{arXiv:math.QA/9808003}.

\bibitem{OO1}
Ohyama Y., Okumura S., A coalescent diagram of the {P}ainlev\'e equations from
 the viewpoint of isomonodromic deformations, \href{https://doi.org/10.1088/0305-4470/39/39/S08}{\textit{J.~Phys.~A}} \textbf{39}
 (2006), 12129--12151.

\bibitem{OO2}
Ohyama Y., Okumura S., R. {F}uchs' problem of the {P}ainlev\'e equations from
 the first to the fifth, in Algebraic and {G}eometric {A}spects of
 {I}ntegrable {S}ystems and {R}andom {M}atrices, \textit{Contemp. Math.}, Vol.
 593, \href{https://doi.org/10.1090/conm/593/11876}{American Mathematical Society}, Providence, RI, 2013, 163--178,
 \href{https://arxiv.org/abs/math.CA/0512243}{arXiv:math.CA/0512243}.

\bibitem{O1}
Okamoto K., Polynomial {H}amiltonians associated with {P}ainlev\'e
 equations.~{I}, \href{https://doi.org/10.3792/pjaa.56.264}{\textit{Proc. Japan Acad. Ser.~A Math. Sci.}} \textbf{56}
 (1980), 264--268.

\bibitem{O}
Okamoto K., Studies on the {P}ainlev\'e equations.~{I}. {S}ixth {P}ainlev\'e
 equation {$P_{{\rm VI}}$}, \href{https://doi.org/10.1007/BF01762370}{\textit{Ann. Mat. Pura Appl.~(4)}} \textbf{146}
 (1987), 337--381.

\bibitem{O2}
Okamoto K., Studies on the {P}ainlev\'e equations.~{II}. {F}ifth {P}ainlev\'e
 equation {$P_{\rm V}$}, \href{https://doi.org/10.4099/math1924.13.47}{\textit{Japan.~J.~Math.~(N.S.)}} \textbf{13} (1987),
 47--76.

\bibitem{Pa-Ra}
Paul E., Ramis J.-P., Dynamics on wild character varieties, \href{https://doi.org/10.3842/SIGMA.2015.068}{\textit{SIGMA}}
 \textbf{11} (2015), 068, 21~pages, \href{https://arxiv.org/abs/1508.03122}{arXiv:1508.03122}.

\bibitem{PR23}
Paul E., Ramis J.-P., Dynamics of the fifth {P}ainlev\'e foliation, \href{https://arxiv.org/abs/2301.08597}{arXiv:2301.08597}.

\bibitem{vdP-Sa}
van~der Put M., Saito M.H., Moduli spaces for linear differential equations and
 the {P}ainlev\'e equations, \href{https://doi.org/10.5802/aif.2502}{\textit{Ann. Inst. Fourier (Grenoble)}}
 \textbf{59} (2009), 2611--2667, \href{https://arxiv.org/abs/0902.1702}{arXiv:0902.1702}.

\bibitem{vdP-S}
van~der Put M., Singer M.F., Galois theory of linear differential equations,
 \textit{Grundlehren Math. Wiss.}, Vol.~328, \href{https://doi.org/10.1007/978-3-642-55750-7}{Springer}, Berlin, 2003.

\bibitem{vdP-T}
van~der Put M., Top J., A {R}iemann--{H}ilbert approach to {P}ainlev\'e~{IV},
 \href{https://doi.org/10.1080/14029251.2013.862442}{\textit{J.~Nonlinear Math. Phys.}} \textbf{20} (2013), 165--177,
 \href{https://arxiv.org/abs/1207.4335}{arXiv:1207.4335}.

\bibitem{vdP-T3}
van~der Put M., Top J., Geometric aspects of the {P}ainlev\'e equations {${\rm
 PIII}(\rm D_6)$} and {${\rm PIII}(\rm D_7)$}, \href{https://doi.org/10.3842/SIGMA.2014.050}{\textit{SIGMA}} \textbf{10}
 (2014), 050, 24~pages, \href{https://arxiv.org/abs/1207.4023}{arXiv:1207.4023}.

\bibitem{S-H-C}
Sen A., Hone A.N.W., Clarkson P.A., On the {L}ax pairs of the symmetric
 {P}ainlev\'e equations, \href{https://doi.org/10.1111/j.1467-9590.2006.00356.x}{\textit{Stud. Appl. Math.}} \textbf{117} (2006),
 299--319.

\end{thebibliography}
\end{document}